\documentclass[final]{siamart0216}
\usepackage{amsfonts}

\usepackage{bm}
\usepackage{amsfonts,amsmath,amssymb}
\usepackage{mathrsfs,mathtools,stmaryrd,wasysym}
\usepackage{enumerate}
\usepackage{enumitem}
\usepackage{esint}
\usepackage{graphicx,psfrag}
\usepackage{hyperref}
\DeclareMathAlphabet{\mathpzc}{OT1}{pzc}{m}{it}

\newcommand{\FF}[1]{{\color{black}#1}}

\newcommand{\EO}[1]{{\color{black}#1}}

\usepackage{lineno}
\usepackage{pgf,tikz,pgfplots}
\usepackage{pstricks-add}

\DeclareFontEncoding{FMS}{}{}
\DeclareFontSubstitution{FMS}{futm}{m}{n}
\DeclareFontEncoding{FMX}{}{}
\DeclareFontSubstitution{FMX}{futm}{m}{n}
\DeclareSymbolFont{fouriersymbols}{FMS}{futm}{m}{n}
\DeclareSymbolFont{fourierlargesymbols}{FMX}{futm}{m}{n}
\DeclareMathDelimiter{\VERT}{\mathord}{fouriersymbols}{152}{fourierlargesymbols}{147}

\newsiamremark{remark}{Remark}

\newcommand{\TheTitle}{Error estimates for optimal control problems \EO{involving} the Stokes system \EO{and} Dirac measures}
\newcommand{\ShortTitle}{PDE--Constrained Optimization with Dirac Measures}
\newcommand{\TheAuthors}{F. Fuica, E. Ot\'arola, D. Quero}

\headers{\ShortTitle}{\TheAuthors}

\title{{\TheTitle}\thanks{\FF{FF is supported by UTFSM through Beca de Mantenci\'on}. EO is partially supported by CONICYT through FONDECYT project 11180193. DQ is partially supported by \FF{UTFSM} through Programa de Incentivos a la Iniciaci\'on Cient\'ifica (PIIC).}}

\author{Francisco Fuica\thanks{Departamento de Matem\'atica, Universidad T\'ecnica Federico Santa Mar\'ia, Valpara\'iso, Chile.
(\email{francisco.fuica@sansano.usm.cl}).}
\and
Enrique Ot\'arola\thanks{Departamento de Matem\'atica, Universidad T\'ecnica Federico Santa Mar\'ia, Valpara\'iso, Chile.
(\email{enrique.otarola@usm.cl}, \url{http://eotarola.mat.utfsm.cl/}).}
\and
Daniel Quero\thanks{Departamento de Matem\'atica, Universidad T\'ecnica Federico Santa Mar\'ia, Valpara\'iso, Chile.
(\email{daniel.quero@alumnos.usm.cl}).}}

\ifpdf
\hypersetup{
  pdftitle={\TheTitle},
  pdfauthor={\TheAuthors}
}
\fi

\date{Draft version of \today.}

\begin{document}

\maketitle

\begin{abstract}
The aim of this work is to derive a priori error estimates for \EO{finite element discretizations of} control--constrained optimal control problems that involve the Stokes system and Dirac measures. The first problem entails the minimization of a cost functional that involves point evaluations of the velocity field that solves the state equations. This leads to an adjoint problem with a linear combination of Dirac measures as a forcing term and whose solution exhibits reduced regularity properties. The second problem involves a control variable that corresponds to the amplitude of forces modeled as point sources. This leads to a solution of the state equations with  reduced regularity properties. For each problem, we propose a \EO{finite element} solution technique and derive \EO{a priori} error estimates. Finally, we present numerical experiments, in two and three dimensions, \EO{that illustrate our theoretical developments.}
\end{abstract}



\section{Introduction}

Let $d \in \{2,3\}$ and  $\Omega\subset\mathbb{R}^{d}$ be an open and bounded polytopal domain with Lipschitz boundary $\partial\Omega$. The purpose of this work is the study of a priori error estimates for finite element \EO{discretizations that approximate the solution to} optimal control problems involving the Stokes equations and Dirac measures; control--constraints are also considered. We consider two illustrative examples, which we proceed to describe in what follows.

\subsection{Optimization with Point Observations}\label{sec:point_observations} Let $\mathcal{Z} \neq\emptyset$ be a finite ordered subset of $\Omega$ with cardinality $\#\mathcal{Z} = m$. Given a set of desired states $\{{\bf{y}}_t\}_{t\in {\mathcal{Z}}}\subset \mathbb{R}^d$, a regularization parameter $\lambda>0$, and the cost functional 
\begin{equation}\label{def:cost_func}
J(\mathbf{y},\mathbf{u}):=\frac{1}{2}\sum_{t\in {\mathcal{Z}}}|\mathbf{y}(t)-\mathbf{y}_t|^2+\frac{\lambda}{2}\|\mathbf{u}\|_{{\mathbf{L}}^2(\Omega)}^2,
\end{equation} 
we are interested in finding $\min J(\mathbf{y},\mathbf{u})$ subject to the Stokes system
\begin{equation}\label{def:state_eq}
\left\{
\begin{array}{rcll}
-\Delta \mathbf{y} + \nabla p & = & \mathbf{u} & \text{ in } \quad \Omega, \\
\text{div }\mathbf{y} & = & 0 & \text{ in } \quad \Omega, \\
\mathbf{y} & = & \mathbf{0} & \text{ on } \quad \partial\Omega,
\end{array}
\right.
\end{equation}
and the control constraints
\begin{equation}\label{def:box_constraints}
\mathbf{u}\in \mathbb{U}_{ad},\quad \mathbb{U}_{ad}:=\{\mathbf{v} \in \mathbf{L}^2(\Omega):  \mathbf{a} \leq \mathbf{v} \leq \mathbf{b} \text{ a.e. in } \Omega \},
\end{equation}
with $\mathbf{a},\mathbf{b} \in \mathbb{R}^{d}$ satisfying \EO{$-\boldsymbol{\infty}<\mathbf{a} < \mathbf{b}<\boldsymbol{\infty}$}. We immediately comment that, throughout this work, vector inequalities must be understood componentwise and that $|\cdot|$ will denote the Euclidean norm in $\mathbb{R}^{d}$.

In contrast to standard PDE--constrained optimization problems, the cost functional \eqref{def:cost_func} involves point evaluations of the state velocity field. This leads to a subtle formulation of the adjoint problem:
\begin{equation}
\label{def:adjoint_eq}
\left\{
\begin{array}{rcll}
-\Delta \mathbf{z} - \nabla r  & = &  \sum_{t \in \mathcal{Z}} (\mathbf{y}\FF{(t)} - \mathbf{y}_t ) \delta_t & \text{ in } \Omega, \\
\text{div }\mathbf{z}  & = & 0  & \text{ in }  \Omega, \\
\mathbf{z}  &= & \mathbf{0}  & \text{ on }  \partial\Omega.
\end{array}
\right.
\end{equation}
Here, $\delta_t$ corresponds to the Dirac delta supported at the interior point $t$.

The optimal control problem \eqref{def:cost_func}--\eqref{def:box_constraints} finds relevance in some applications where the observations are carried out at specific locations of the domain. For instance, in the active control of sound \cite{PANelson,MR2086168} and in the active control of vibrations \cite{1996ix,MR2525606}; see also \cite{MR3523574,MR3449612} for other applications. We immediately comment that, since the domain $\Omega$ is Lipschitz and $\mathbb{U}_{ad} \subset \mathbf{L}^{\infty}(\Omega)$, the point observations \EO{$ \mathbf{y}(t)$, with $t \in \mathcal{Z}$,} are well--defined; see \cite[Theorem 2.9]{Brown_Shen} and \cite[Lemma 12]{demlar13}. \EO{These point observations} tend to enforce the state velocity field $\mathbf{y}$ to have the fixed vector value $\mathbf{y}_{t}$ at the point \EO{$t \in \mathcal{Z}$}. Consequently, problem \eqref{def:cost_func}--\eqref{def:box_constraints} can be understood as a penalty version of a PDE--constrained optimization problem where the velocity field that solves the state equation is constrained at a collection of points. 

There are several works that provide a priori error estimates for the optimal control problem \eqref{def:cost_func}--\eqref{def:box_constraints} \EO{but} when the state equations are a Poisson problem. Under the fact that the associated adjoint variable belongs to $W_0^{1,r}(\Omega)$, for $r \in (1,d/(d-1))$, the authors of \cite{MR3449612} obtain, for $d \in \{ 2,3\}$, a priori and a posteriori error estimates for the so-called variational discretization scheme when applied for discretizing the underlying optimal control problem; the state and adjoint equations are discretized on the basis of standard piecewise linear finite element functions. The derived a priori error estimates for the control, the state, and adjoint state variables are optimal in terms of regularity; see \cite[Theorem 3.2]{MR3449612}. Later, the authors of \cite{MR3523574} propose and analyze a fully discrete scheme that discretizes the optimal state, adjoint, and control variables with piecewise linear functions. For this scheme, the authors provide error estimates when $d=2$ \cite[Theorem 5.1]{MR3523574}; the control and the state are discretized using meshes of size $\mathcal{O}(h^2)$ and $\mathcal{O}(h)$, respectively. The authors of \cite{MR3523574} also analyze the variational discretization scheme and derive an optimal error estimate, in terms of regularity, for the control variable \cite[Theorem 5.2]{MR3523574}. In \cite{MR3800041}, the authors invoke the theory of Muckenhoupt weights and Muckenhoupt--weighted Sobolev spaces to provide error estimates for a numerical scheme that discretize\FF{s} the control variable with piecewise constant functions;  the state and adjoint equations are discretized with piecewise linear finite elements. To be precise, the authors derive a priori error estimates for the error approximation of the optimal control variable when $d \in \{2,3\}$; the one for $d=2$ being nearly--optimal in terms of approximation \cite[Theorem 4.3]{MR3800041}. However, the estimate for $d=3$ is suboptimal in terms of approximation; it behaves as $\mathcal{O}(h^{1/2}|\log h|)$. This has been recently improved in \cite[Theorem 6.6]{2018arXiv180202918B}.


\subsection{Optimization with Singular Sources}
\label{sec:singular_source} 
Let $\mathcal{D}\neq\emptyset$ be a finite ordered subset of $\Omega$ with cardinality $\#\mathcal{D} = l$. Given a desired state ${\bf{y}}_\Omega\in \mathbf{L}^2(\Omega)$, a regularization parameter $\lambda>0$, and the cost functional 
\begin{equation}\label{def:cost_func_2}
\mathfrak{J}(\mathbf{y},\mathcal{U}):=\frac{1}{2}\|\mathbf{y}-\mathbf{y}_\Omega\|_{\mathbf{L}^2(\Omega)}^{2}+\frac{\lambda}{2}\sum_{t\in\mathcal{D}}|\mathbf{u}_t|^2, \quad \mathcal{U}= (\mathbf{u}_{1},\ldots,\mathbf{u}_{l}),
\end{equation} 
the problem under consideration reads as follows: Find $\min \mathfrak{J}(\mathbf{y},\mathcal{U})$ subject to 
\begin{equation}\label{def:state_eq_2}
-\Delta \mathbf{y} + \nabla p  =  \sum_{t\in\mathcal{D}}\mathbf{u}_t \delta_t \text{ in } \Omega, 
\quad
\text{div }\mathbf{y}  =  0  \text{ in } \Omega, 
\quad
\mathbf{y}  =  \mathbf{0}  \text{ on } \partial\Omega,
\end{equation}
and the control constraints $\mathcal{U}=(\mathbf{u}_1,\ldots,\mathbf{u}_{l})\in \mathfrak{U}_{ad}$, where
\begin{equation}\label{def:box_constraints_2}
\mathfrak{U}_{ad}:=\{\mathcal{V}=(\mathbf{v}_1,\ldots,\mathbf{v}_{l}) \in [\mathbb{R}^d]^{l}:  \mathbf{a}_t \leq \mathbf{v}_t \leq \mathbf{b}_t \text{ for all } t\in \mathcal{D} \},
\end{equation}
with $\mathbf{a}_t,\mathbf{b}_t \in \mathbb{R}^{d}$ satisfying $\mathbf{a}_t < \mathbf{b}_t$  for all $t\in\mathcal{D}$. 

The optimization problem \eqref{def:cost_func_2}--\eqref{def:box_constraints_2} is of relevance in applications where it can be specified a control at finitely many prespecified points. We observe that, in view of the particular structure of the control variable $\mathcal{U}$, the state equation \eqref{def:state_eq_2} corresponds to a Stokes system that has a linear combination of Dirac measures on the right--hand side of the \emph{momentum equation}. 

There are a few works that consider the numerical approximation of problem \eqref{def:cost_func_2}--\eqref{def:box_constraints_2} when the Stokes equations are replaced by a Poisson problem. In \cite{MR3225501}, the authors use the variational discretization concept to derive a priori error estimates. Their technique is based on the fact that the state belongs to $W_0^{1,r}(\Omega)$ for $r \in (1,d/(d-1))$. An approach involving weighted estimates has also been considered in \cite{MR3800041}, where the authors obtain the following rates of convergence for the error approximation of the control variable: $\mathcal{O}(h^{2-\epsilon})$ in two dimensions and $\mathcal{O}(h^{1-\epsilon})$ in three dimensions, where $\epsilon > 0$.

Since the Stokes system with a linear combination of Dirac measures in the \emph{momentum equation} appears as the adjoint system \eqref{def:adjoint_eq} of problem \eqref{def:cost_func}--\eqref{def:box_constraints} and as the state equation of problem \eqref{def:cost_func_2}--\eqref{def:box_constraints_2}, it is of importance to understand the regularity properties of the involved solution and the development of numerical techniques to approximate it. The main difficulty in the study of the aforementioned problem is that it does not admit a solution in the classical Hilbert space $\mathbf{H}^{1}_{0}(\Omega) \times L^{2}(\Omega)/\mathbb{R}$. In spite of this fact, there are a few works that consider the numerical approximation of this system. In \cite{LACOUTURE2015187} the author presents a numerical method in two and three--dimensional bounded domains. \EO{However,} convergence properties are not investigated. Later, the authors of \cite{MR3854357} derive quasi--optimal local convergence results in $\mathbf{H}^1 \times L^2$; the error is analyzed on a subdomain which does not contain the singularity of the solution. The authors operate in two dimensions and consider the mini element and Taylor--Hood schemes. On the other hand, on the basis of the fact that the solution to the Stokes system with singular sources can be seen as an element of a weighted space, a priori and a posteriori error estimates, for classical low--order inf--sup stable finite element approximations, have been developed in \cite{2019arXiv190500476D} and \cite{Allendes_et_al2017}, respectively. \EO{The a priori error estimates derived in \cite[Corollary 5.4]{2019arXiv190500476D} have been recently improved in \cite{Otarola_Salgado2020}.}

In spite of these advances and to the best of our knowledge, this is the first work that provides a priori error estimates for the optimal control problems \eqref{def:cost_func}--\eqref{def:box_constraints} and 
\eqref{def:cost_func_2}--\eqref{def:box_constraints_2}. We discretize the adjoint and state equations with classical low--order inf--sup stable finite element \EO{schemes} and the control variable with piecewise constant functions. \EO{We} derive error estimates for problem \eqref{def:cost_func}--\eqref{def:box_constraints}, \EO{in two and three--dimensions, and problem \eqref{def:cost_func_2}--\eqref{def:box_constraints_2}, in two--dimensions. In particular, in two--dimensions, and for both problems, the derived error estimates for the discretization of the optimal control variable are nearly--optimal in terms of approximation. In three--dimensions, and for problem \eqref{def:cost_func}--\eqref{def:box_constraints}, the obtained error estimate for the control variable is suboptimal in terms of approximation. However, numerical experiments seem to indicate that such an estimate is attained in practice.}

The rest of the paper is organized as follows. In Section \ref{sec:notation_and_prel} we introduce the notation and functional framework we shall work with. We also briefly review, in Section \ref{sec:Stokes_dirac}, the well--posedness of the Stokes system with singular sources. Section \ref{sec:pointwise_tracking} contains the numerical analysis for problem \eqref{def:cost_func}--\eqref{def:box_constraints}. In Section \ref{sec:singular_sources_control_probl} we derive error estimates for the optimal control problem \eqref{def:cost_func_2}--\eqref{def:box_constraints_2}. We conclude, in Section \ref{sec:numerical_ex}, with a series of numerical examples that illustrate the developed theory.


\section{Notation and Preliminaries}\label{sec:notation_and_prel}
Let us fix the notation and conventions in which we will operate. 

\subsection{Notation}\label{sec:notation}
Throughout this work $d\in\{2,3\}$ and $\Omega\subset\mathbb{R}^d$ is an open, bounded, and convex polytopal domain. If $\mathcal{X}$ and $\mathcal{Y}$ are normed vector spaces, we write $\mathcal{X}\hookrightarrow\mathcal{Y}$ to denote that $\mathcal{X}$ is continuously embedded in $\mathcal{Y}$. We denote by $\mathcal{X}'$ and $\|\cdot\|_\mathcal{X}$ the dual and the norm of $\mathcal{X}$, respectively. Given a Lebesgue measurable subset $A$ of $\mathbb{R}^d$, we denote by $|A|$ its Lebesgue measure.

To denote vector-valued functions we shall use lower--case bold letters, whereas to denote function spaces we shall use upper-case bold letters. For a bounded domain $G \subset \mathbb{R}^d$, if $X(G)$ corresponds to a function space over $G$, we shall denote $\mathbf{X}(G)=[X(G)]^{d}$. In particular, we denote $\mathbf{L}^2(G)=[L^2(G)]^d$, which is equipped with the following inner product and norm, respectively:
\begin{equation*}
(\mathbf{w},\mathbf{v})_{\mathbf{L}^2(G)}=\int_G \mathbf{w}\cdot \mathbf{v}, 
\qquad
\|\mathbf{v}\|_{\mathbf{L}^2(G)}=(\mathbf{v},\mathbf{v})_{\mathbf{L}^2(G)}^{\frac{1}{2}}\qquad \forall\: \mathbf{w},\mathbf{v}\in \mathbf{L}^2(G).
\end{equation*}

Finally, the relation $a \lesssim b$ indicates that $a \leq C b$, with a positive constant that does not depend on $a$, $b$ nor the discretization parameter. The value of $C$ might change at each occurrence.


\subsection{Weighted Sobolev Spaces}\label{sec:wse}

We start this section with a notion which will be fundamental for further discussions, that of a weight. A weight is a nonnegative locally integrable function on $\mathbb{R}^d$ that takes values in $(0,\infty)$ almost everywhere. We will be particularly interested in the weights belonging to the so-called Muckenhoupt class $A_2(\mathbb{R}^d)$ \cite{Javier2001,Fabes_et_al1982,MR0293384,Turesson2000}.

\begin{definition}[Muckenhoupt class $A_2(\mathbb{R}^d)$]
\label{def:muckenhoupt_a2}
Let $\omega$ be a weight. We say that $\omega \in A_{2}(\mathbb{R}^d)$, or that $\omega$ is an $A_{2}(\mathbb{R}^d)$--weight, if there exists a positive constant $C_{\omega}$ such that
\begin{equation*}
C_{\omega} = \sup_{B}\left(\frac{1}{|B|}\int_B\omega \right)\left(\frac{1}{|B|}\int_B\omega^{-1} \right) < \infty,
\end{equation*}
where the supremum is taken over all balls $B$ in $\mathbb{R}^{d}$.
\end{definition}

Let $\omega \in A_2(\mathbb{R}^d)$ and $G\subset \mathbb{R}^{d}$ be an open and bounded domain. We define the weighted Lebesgue space $L^2(\omega,\EO{G})$ as the space of square--integrable functions with respect to the measure $\omega \textrm{d}x$.
We also define the weighted Sobolev space
$
H^1(\omega,G):=\{v\in L^2(\omega,G) : |\nabla v| \in L^2(\omega,G) \},
$
which we equip with the norm
\begin{equation*}
\|v\|_{H^1(\omega,G)}:= \left(\|v\|_{L^2(\omega,G)}^2 + \|\nabla v\|_{L^2(\omega,G)}^2\right)^\frac{1}{2}.
\end{equation*}
It is remarkable that most of the properties of classical Sobolev spaces have a weighted counterpart. This is not because of the specific form of the weight but rather due to the fact that the weight $\omega \in A_2(\mathbb{R}^d)$. In particular, $L^{2}(\omega,G)$ and $H^1(\omega,G)$ are Hilbert spaces and  $C^\infty(\FF{G})$ is dense in $H^1(\omega,G)$; see, for instance, \cite[Proposition 2.1.2, Corollary 2.1.6]{Turesson2000} and \cite[Theorem 1]{Goldshtein_Ukhlov_2009}. Define $H^{1}_{0}(\omega,G)$ as the closure of $C_{0}^{\infty}(\FF{G})$ in $H^1(\omega,G)$. In view of a weighted Poincar\'e inequality, that follows from \cite[Theorem 1.3]{Fabes_et_al1982} for $\FF{G}$ being a ball and \cite{MR1140667,MR1415588} for more general domains, we conclude that, in $H^{1}_{0}(\omega,G)$, the seminorm $\|\nabla v\|_{L^2(\omega,G)}$ is equivalent to the norm $\|v\|_{H^1(\omega,G)}$.

Finally, we define the vector space $\mathbf{H}^{1}_{0}(\omega,G):=[H^{1}_{0}(\omega,G)]^d$, which, in view of the aforementioned Poincar\'e inequality, we equip with the norm
\begin{equation*}
\|\nabla \mathbf{v} \|_{\mathbf{L}^{2}(\omega,G)} = \left(\sum_{i=1}^{d}\|\nabla v_{i}\|_{L^{2}(\omega,G)}^2 \right)^{\frac{1}{2}}.
\end{equation*}

\subsection{The Stokes System with Dirac Sources}\label{sec:Stokes_dirac}
In this section we review well--posedness results in weighted spaces for the Stokes system with a linear combination of Dirac measures as a forcing term in the momentum equation. We also comment on the finite element approximation of such a problem. 
The review is motivated since the adjoint equation \eqref{def:adjoint_eq} for the pointwise tracking optimal control problem of Section \ref{sec:point_observations} and the state equation \eqref{def:state_eq_2} for the optimization with singular sources problem of Section \ref{sec:singular_source} are particular instances of the aforementioned singular setting.

Let $\mathcal{E}$ be a finite ordered subset of $\Omega$ with cardinality $\#\mathcal{E}$. Given $\{\mathbf{F}_t\}_{t\in\mathcal{E}}\subset\mathbb{R}^d$, we consider the following boundary value problem: Find $(\boldsymbol\Phi,\zeta)$ such that 
\begin{equation}\label{eq:stokes_delta}
\left\{
\begin{array}{rcll}
-\Delta \boldsymbol\Phi + \nabla \zeta & = & \sum_{t\in\mathcal{E}}\mathbf{F}_t\delta_{t} & \text{ in } \quad \Omega, \\
\text{div}\: \boldsymbol\Phi & = & 0 & \text{ in } \quad \Omega, \\
\boldsymbol\Phi & = & \mathbf{0} & \text{ on } \quad \partial\Omega,
\end{array}
\right.
\end{equation}
where $\delta_{t}$ denotes the Dirac delta supported at ${t}\in\Omega$. Let us assume that $\Omega = \mathbb{R}^d$. If this is the case, the results of \cite[Section IV.2]{Gal11} yield the following asymptotic behavior near the points $t \in \mathcal{E}$:
\begin{equation*}
|\nabla \boldsymbol\Phi(x)| \approx |x-t|^{1-d}, \quad  |\zeta(x)| \approx |x-t|^{1-d}.
\end{equation*}
Consequently, $(\boldsymbol\Phi,\zeta)\not\in \mathbf{H}^1(\Omega)\times L^2(\Omega)/\mathbb{R}$. However, if $B(t^*,r)$ denotes a ball of radius $r$ and center $t^*$ with $t^* \in \mathcal{E}$ and $r>0$  such that $B(t^*,r) \cap \mathcal{E} = \{ t^*\}$, then
\[
\int_{B(t^*,r)} |x-t^*|^{\alpha} |\nabla \boldsymbol\Phi(x)|^2 {\rm d} x < \infty,
\quad
\int_{B(t^*,r)} |x-t^*|^{\alpha} |\zeta(x)|^2 {\rm d} x < \infty,
\]
for $\alpha \in (d-2,\infty)$. This heuristic argument suggests to study problem \eqref{eq:stokes_delta} on weighted spaces.

\subsubsection{Weak Formulation}
\label{subsub:weak_formulation}

Define
\[
d_{\mathcal{E}} := \left\{\begin{array}{ll}
\mathrm{dist}(\mathcal{E},\partial \Omega), & \mbox{if } \# \mathcal{E}=1,
\\
\min \left \{ \mathrm{dist}(\mathcal{E},\partial \Omega), \min\{|t-t'|: t,t' \in \mathcal{E}, \ t\neq t' \} \right \}, & \mbox{otherwise}.
\end{array}\right.
\]
Since $\mathcal{E} \subset \Omega$ and $\mathcal{E}$ is finite, we thus have that $d_{\mathcal{E}}>0$. With this notation, we define the weight $\rho$ as follows: if $ \# \mathcal{E}=1$, then
\begin{equation}\label{def:weight_rho}
\rho(x) = \mathsf{d}_t^\alpha(x),
\end{equation}
otherwise
\begin{equation}\label{def:weight_rho_complete}
\rho(x)=
\begin{cases}
\mathsf{d}_t^\alpha(x),\:&\exists t\in \mathcal{E}:\mathsf{d}_t(x) < \frac{d_\mathcal{E}}{2},\\
1, &\mathsf{d}_t(x)\geq \frac{d_\mathcal{E}}{2} \: \forall \: t\in \mathcal{E},
\end{cases}
\end{equation}
where $\mathsf{d}_t^\alpha(x) := |x - t|^\alpha$ and $\alpha \in (-d,d)$. Since $\alpha \in (-d,d)$,
owing to \cite[Theorem 6]{ACDT2014} and \cite[Lemma 2.3 (v)]{MR1601373}, we have that the function $\rho$ belongs to the Muckenhoupt class $A_{2}(\mathbb{R}^d)$. Define the spaces
\begin{align*}
\mathbb{X}:&= \mathbf{H}_0^1(\Omega) + \mathbf{H}_0^1(\rho,\Omega)+\mathbf{H}_0^1(\rho^{-1},\Omega),\\
\mathbb{Y}:&=L^2(\Omega)/\mathbb{R}+ L^2(\rho,\Omega)/\mathbb{R}+ L^2(\rho^{-1},\Omega)/\mathbb{R}.
\end{align*}
Finally, we define the bilinear forms
\begin{equation*}
a: \mathbb{X}\times \mathbb{X} \rightarrow \mathbb{R}, 
\quad 
a(\mathbf{w},\mathbf{v}):=\int_\Omega \nabla \mathbf{w} : \nabla \mathbf{v}, 
\end{equation*}
and
\begin{equation*}
b:\mathbb{X}\times \mathbb{Y}, 
\quad 
b(\mathbf{v},q):=-\int_\Omega q \:\text{div}\: \mathbf{v}.
\end{equation*}

With all these ingredients at hand, we present the following weak formulation of problem \eqref{eq:stokes_delta} \cite[Section 3]{Allendes_et_al2017}: Find $(\boldsymbol\Phi,\zeta)\in \mathbf{H}_0^1(\rho,\Omega)\times L^2(\rho,\Omega)/\mathbb{R}$ such that 
\begin{equation}\label{eq:weak_stokes_delta}
\left\{\!\!
\begin{array}{rcll}
a(\boldsymbol\Phi,\mathbf{v})+b(\mathbf{v},\zeta) & = & \sum_{t\in\mathcal{E}}\langle \mathbf{F}_t\delta_{t},\mathbf{v}\rangle  &\quad \forall  \mathbf{v} \in  \mathbf{H}_0^1(\rho^{-1},\Omega), \\
b(\boldsymbol\Phi,q) & = & 0 &\quad \forall q \in L^2(\rho^{-1},\Omega)/\mathbb{R},
\end{array}
\right.
\hspace{-0.4cm}
\end{equation}
where $\langle \cdot,\cdot \rangle$ denotes the duality pairing between $\mathbf{H}_0^1(\rho^{-1},\Omega)'$ and $\mathbf{H}_0^1(\rho^{-1},\Omega)$. We immediately mention that in order to guarantee that $\mathbf{F}_t\delta_{t} \in \mathbf{H}_0^1(\rho^{-1},\Omega)'$, and thus that $\langle \mathbf{F}_t\delta_{t},\mathbf{v} \rangle$ is well defined for $\mathbf{v} \in \mathbf{H}_0^1(\rho^{-1},\Omega)$, the parameter $\alpha$ should be restricted to $(d-2,d)$; see  
\cite[Lemma 7.1.3]{Kozlov_et_al1997}, \cite[Remark 21.19]{MR2305115}, and \EO{\cite[Proposition 5.2]{2019arXiv190500476D}} for details.

It can be proved that problem \eqref{eq:weak_stokes_delta} admits a unique solution; see \cite[Theorem 14]{Otarola_Salgado2017}. Moreover, the following estimate can be obtained \cite[Theorem 14]{Otarola_Salgado2017}:
\begin{equation}\label{eq:a_priori_estimates}
\|\nabla \boldsymbol\Phi\|_{\mathbf{L}^2(\rho,\Omega)}+\|\zeta\|_{L^2(\rho,\Omega)/\mathbb{R}}
\lesssim
\sum_{t\in\mathcal{E}}|\mathbf{F}_t| \|\delta_{t}\|_{H_0^1(\rho^{-1},\Omega)'}.
\end{equation}

We conclude this section with the following embedding result.

\begin{theorem}[$\mathbf{H}_0^1(\rho,\Omega) \hookrightarrow \mathbf{L}^2(\Omega)$]
\label{thm:weighted_poincare}
 If $\alpha \in (d - 2,2)$, then $\mathbf{H}_0^1(\rho,\Omega) \hookrightarrow \mathbf{L}^2(\Omega)$. Moreover, the following weighted Poincar\'e inequality holds
\begin{equation*}
\|\mathbf{v}\|_{\mathbf{L}^{2}(\Omega)} \lesssim \|\nabla \mathbf{v}\|_{\mathbf{L}^{2}(\rho,\Omega)} \quad \forall  \mathbf{v} \in \mathbf{H}_{0}^{1}(\rho,\Omega),
\end{equation*}
where the hidden constant depends only on $\Omega$ and $d_{\mathcal{E}}$.
\end{theorem}
{\it Proof.}
The proof follows from \cite[Lemmas 1 and 2]{Allendes_et_al2017_2}.
\qed


\subsubsection{Finite Element Approximation and Error Estimates}
\label{sec:approx_results}
We start the discussion by introducing some standard finite element notation \cite{brenner,CiarletBook,Guermond-Ern}.  We denote by $\mathscr{T}_h = \{ T\}$ a conforming partition, or mesh, of $\bar{\Omega}$ into closed simplices $T$ with size $h_T = \text{diam}(T)$. Define $h:=\max_{ T \in \mathscr{T}_h} h_T$. \EO{We denote by} $\mathbb{T} = \{\mathscr{T}_h \}_{h>0}$ a collection of conforming and quasi--uniform meshes. \EO{For $T \in \mathscr{T}_h$, we define the \emph{patch} $\mathcal{N}_T$ associated with the element $T$ as
\begin{equation}\label{def:patch}
\mathcal{N}_T:= \left \{ T^{\prime}\in\mathscr{T}: {T}\cap {T^\prime}\neq\emptyset \right \}.
\end{equation}
In an abuse of notation, in what follows, by $\mathcal{N}_T$ we will indistinctively denote either this set or the union of the triangles that comprise it.} 

Given a mesh $\mathscr{T}_h \in \mathbb{T}$, we denote by $\mathbf{V}_h$ and $Q_h$ the finite element spaces that approximate the velocity field and the pressure, respectively, constructed over $\mathscr{T}_h$. In this work we will consider the following popular finite element discretizations: 
\begin{itemize}
\item[(a)] The mini element \cite[Section 4.2.4]{Guermond-Ern}: 
\begin{align}\label{def:discrete_spaces_P1B}
\begin{split}
Q_h &=\left\{ q_h \in C(\bar{\Omega})\, : \, q_h|^{}_T \in \mathbb{P}_{1}(T) \ \forall \: T \in\mathscr{T}_h \right\}  \cap L^2(\Omega)/\mathbb{R}, \\
\mathbf{V}_h &=\left\{ \mathbf{v}_h \in \mathbf{C}(\bar{\Omega})\, : \, \mathbf{v}_h|^{}_T \in [ \mathbb{P}_1(T) \oplus \mathbb{B}(T) ]^d \ \forall \: T\in\mathscr{T}_h \right\}\cap \mathbf{H}_0^1(\Omega),
\end{split}
\end{align}
where $\mathbb{B}(T)$ denotes the space spanned by local bubble functions.
\item[(b)] The classical Taylor--Hood elements \cite[Section 4.2.5]{Guermond-Ern}:
\begin{align}\label{def:discrete_spaces_TH}
\begin{split}
Q_h &=\left\{ q_h \in C(\bar{\Omega})\, : \, q_h|^{}_T \in \mathbb{P}_{1}(T) \ \forall \: T \in\mathscr{T}_h \right\}  \cap L^2(\Omega)/\mathbb{R}, \\
\mathbf{V}_h &=\left\{ \mathbf{v}_h \in \mathbf{C}(\bar{\Omega})\, : \, \mathbf{v}_h|^{}_T \in [\mathbb{P}_2(T)]^d \ \forall \: T\in\mathscr{T}_h \right\}\cap \mathbf{H}_0^1(\Omega).
\end{split}
\end{align}
\end{itemize}

We now present an error estimate.

\begin{lemma}[Error estimate for Stokes system with Dirac sources]
\label{lemma:first_estimate}
Let $\Omega$ be convex. Let $(\boldsymbol\Phi,\zeta)\in\mathbf{H}_0^1(\rho,\Omega)\times L^2(\rho,\Omega)/\mathbb{R}$ be the solution of \eqref{eq:weak_stokes_delta}. Let $(\boldsymbol\Phi_h,\zeta_h)\in \mathbf{V}_h\times Q_h$ be the finite element approximation of $(\boldsymbol\Phi,\zeta)$ on the basis of the discrete spaces \eqref{def:discrete_spaces_P1B} or \eqref{def:discrete_spaces_TH}. \EO{Then, we have the error estimate}
\begin{equation}\label{eq:apriori_error_delta}
\|\boldsymbol\Phi-\boldsymbol\Phi_h\|_{\mathbf{L}^2(\Omega)}
\lesssim
h^{2-d/2}\sum_{t\in\mathcal{E}}|\mathbf{F}_t| \|\delta_{t}\|_{\mathcal{M}(\Omega)},
\end{equation}
\EO{with a hidden constant that does not depend on $(\boldsymbol\Phi,\zeta)$, $(\boldsymbol\Phi_h,\zeta_h)$, nor $h$. Here, $\mathcal{M}(\Omega)$ denotes the space of bounded Radon measures.} 
\end{lemma}
{\it Proof.}
See \EO{\cite[Theorem 7]{Otarola_Salgado2020}}. 
\qed


\section{The Pointwise Tracking Optimal Control Problem}\label{sec:pointwise_tracking}

In this section we analyze a weak version of the optimal control problem \eqref{def:cost_func}--\eqref{def:box_constraints}, which reads \EO{as follows}: Find 
\begin{equation}\label{def:weak_ocp}
\min \{ J(\mathbf{y},\mathbf{u}): (\mathbf{y},\mathbf{u}) \in \mathbf{H}_0^1(\Omega)\times\mathbb{U}_{ad} \},
\end{equation}
subject to 
\begin{equation}\label{eq:weak_pde}
\left\{\begin{array}{rcll}
a(\mathbf{y},\mathbf{v})+b(\mathbf{v},p) & = & (\mathbf{u},\mathbf{v})_{{\mathbf{L}}^2(\Omega)} & \quad \forall  \mathbf{v} \in \mathbf{H}_0^1(\Omega), \\
b(\mathbf{y},q) & = & 0 & \quad \forall  q \in L^2(\Omega)/\mathbb{R}.
\end{array}
\right.
\end{equation}

Standard arguments that rely on the coercivity of $a$ on $\mathbf{H}_0^1(\Omega)$ and \EO{an} inf--sup condition for $b$ yield the existence of a unique solution $(\mathbf{y},p) \in \mathbf{H}_0^1(\Omega) \times L^2(\Omega)/\mathbb{R}$ to problem \eqref{eq:weak_pde}. \EO{In addition,} the pair $(\mathbf{y},p)$ satisfies the following stability estimate:
\begin{equation}\label{eq:standard_estimate_stokes}
\| \nabla \mathbf{y}\|_{\mathbf{L}^2(\Omega)} + \|p\|_{L^{2}(\Omega)} \lesssim \|\mathbf{u}\|_{\mathbf{L}^{2}(\Omega)};
\end{equation}
see, for instance, \cite[Theorem 4.3]{Guermond-Ern}. We also present the following regularity result for the pair $(\mathbf{y},p)$ that solves \eqref{eq:weak_pde}.

\begin{proposition}[regularity]
Let $\Omega \subset \mathbb{R} ^d$ be a convex polytope and $\mathbf{u} \in \mathbf{L}^2(\Omega)$. If $(\mathbf{y},p)$ solves \eqref{eq:weak_pde}, then $\mathbf{y} \in\mathbf{H}^2(\Omega) \cap \mathbf{H}_0^1(\Omega)$, $p \in H^1(\Omega) \cap L^2(\Omega) / \mathbb{R}$, and
\begin{equation}\label{eq:H2regularity}
\|  \mathbf{y}\|_{\mathbf{H}^2(\Omega)} + \|p\|_{H^{1}(\Omega)} \lesssim \|\mathbf{u}\|_{\mathbf{L}^{2}(\Omega)},
\end{equation}
with a hidden constant that is independent of $(\mathbf{y},p)$ and $\mathbf{u}$.
\label{pro:H2regularity}
\end{proposition}
{\it Proof.}
See  \cite[Theorem 2]{MR0404849} and \cite{MR775683} for $d = 2$, and \cite{MR977489} \EO{and \cite{MR1301452}} for $d=3$; see also \cite[Corollary 1.8]{MR2987056}.
\qed

Note that, since the control variable $\mathbf{u}\in \mathbb{U}_{ad}$ and $\Omega$ is convex, the results of Proposition \ref{pro:H2regularity} guarantee that $\mathbf{y}\in \mathbf{H}^2(\Omega)\hookrightarrow \mathbf{C}(\Omega)$. As a consequence, point evaluations of the velocity field $\mathbf{y}$ on the cost functional $J$ are well defined.

\subsection{\EO{Optimality Conditions}}

Since $\lambda>0$ and the underlying control-to-state operator $\mathbf{S}$ is linear and continuous, standard arguments yield the existence of a unique solution $(\bar{\mathbf{y}},\bar{\mathbf{u}})\in\mathbf{H}_0^1(\Omega)\times \mathbb{U}_{ad}$ to the optimal control problem \eqref{def:weak_ocp}--\eqref{eq:weak_pde}.
To present optimality conditions, we introduce the adjoint pair
$(\mathbf{z},r)$ as the unique solution to the following problem: Find $(\mathbf{z},r) \in \mathbf{H}_{0}^{1}(\rho,\Omega) \times L^{2}(\rho,\Omega)/\mathbb{R}$ such that
\begin{equation}
\label{eq:adj_eq}
\left\{\begin{array}{rcll}
a(\mathbf{w},\mathbf{z}) - b(\mathbf{w},r) & = & \displaystyle{\sum_{t\in\mathcal{Z}}}\langle ({\mathbf{y}}\FF{(t)} - \mathbf{y}_{t})\delta_{t},\mathbf{w} \rangle & \quad \forall \mathbf{w} \in \mathbf{H}_0^1(\rho^{-1},\Omega), \\
b(\mathbf{z},s) & = & 0 & \quad \forall   s \in L^2(\rho^{-1},\Omega)/\mathbb{R},
\end{array}
\right.
\end{equation}
where $\mathbf{y}$ solves \eqref{eq:weak_pde}. The weight $\rho$ is defined as in \eqref{def:weight_rho}--\eqref{def:weight_rho_complete} with obvious modifications.
The following first--order sufficient and necessary optimality condition follows from \cite[Theorem 8]{MR4013930}: $(\bar{\mathbf{y}},\bar{\mathbf{u}})$ is optimal for \eqref{def:weak_ocp}--\eqref{eq:weak_pde} if and only if $\bar{\mathbf{u}}$ satisfies
\begin{equation}\label{eq:variational_ineq}
(\bar{\mathbf{z}} + \lambda \bar{\mathbf{u}},\mathbf{u} - \bar{\mathbf{u}})_{\mathbf{L}^{2}(\Omega)} \geq 0  \quad \forall \mathbf{u} \in \mathbb{U}_{ad}.
\end{equation}
Here, $(\bar{\mathbf{z}},\bar{r}) \in \mathbf{H}^{1}_{0}(\rho,\Omega) \times L^{2}(\rho,\Omega)/ \mathbb{R}$ denotes  the optimal adjoint state, which solves \eqref{eq:adj_eq} with $\mathbf{y}$ replaced by $\bar{\mathbf{y}}$. The well--posedness of problem \eqref{eq:adj_eq} follows from the results elaborated in Section \ref{sec:Stokes_dirac}.

We now recall the so-called projection formula for $\bar{\mathbf{u}}$: The optimal control $\bar{\mathbf{u}}$ satisfies \eqref{eq:variational_ineq} if and only if \begin{equation}
\label{eq:projection_formula}
\bar{\mathbf{u}}=\Pi_{[\mathbf{a},\mathbf{b}]}\left(-\lambda^{-1}\bar{\mathbf{z}}\right)  \textrm{ a.e. in } \Omega,
\end{equation}
where $\Pi_{[\mathbf{a},\mathbf{b}]}:\mathbf{L}^1(\Omega)\rightarrow \mathbb{U}_{ad}$ is such that
$
\Pi_{[\mathbf{a},\mathbf{b}]}(\mathbf{v}) := \min\{\mathbf{b},\max\{\mathbf{v},\mathbf{a}\}\}.
$
This projection formula leads to the following regularity result for $\bar{\mathbf{u}}$.
\begin{proposition}(Regularity of $\bar{\mathbf{u}}$)
\label{prop:control_regularity}
If $\bar{\mathbf{u}}$ is optimal for problem \eqref{def:weak_ocp}--\eqref{eq:weak_pde}, then $\bar{\mathbf{u}} \in \mathbf{H}^{1}(\rho,\Omega)$. Moreover, the following estimate holds:
\begin{equation*}
\|\nabla \bar{\mathbf{u}}\|_{\mathbf{L}^2(\rho,\Omega)}
\lesssim
\|\nabla \bar{\mathbf{z}}\|_{\mathbf{L}^2(\rho,\Omega)}.
\end{equation*}
\end{proposition}
{\it Proof.}
Note that, in view of \eqref{eq:projection_formula}, $\bar{\mathbf{u}}$ can be written as
\[\bar{\mathbf{u}} = -\lambda^{-1}\bar{\mathbf{z}} + \max\{\mathbf{a} +\lambda^{-1}\bar{\mathbf{z}},\mathbf{0}\} - \max\{-\lambda^{-1}\bar{\mathbf{z}} - \mathbf{b}, \mathbf{0}\}.\]
The regularity result thus follows directly from \cite[Theorem A.1]{MR567696}. This concludes the proof.
\qed

To summarize, the pair $(\bar{\mathbf{y}},\bar{\mathbf{u}})$ is optimal for the pointwise tracking optimal control problem \eqref{def:weak_ocp}--\eqref{eq:weak_pde} if and only if $(\bar{\mathbf{y}},\bar{p},\bar{\mathbf{z}},\bar{r},\bar{\mathbf{u}})\in \mathbf{H}_0^1(\Omega)\times L^2(\Omega)/\mathbb{R}\times \mathbf{H}_0^1(\rho,\Omega)\times L^2(\rho,\Omega)/\mathbb{R}\times \mathbb{U}_{ad}$ solves the optimality system \eqref{eq:weak_pde}, \eqref{eq:adj_eq}, and \eqref{eq:variational_ineq}.


\subsection{\EO{A Fully Discrete Scheme: Error Estimates}}
\label{sec:fully_discrete_scheme}

\EO{In this section, we design and analyze a fully discrete scheme for the numerical approximation of
problem \eqref{def:weak_ocp}--\eqref{eq:weak_pde}. To begin with,} we define the discrete admissible set
\begin{equation*}
\mathbb{U}_{ad}^h := \mathbf{U}_h \cap \mathbb{U}_{ad},
\end{equation*}
where $\mathbf{U}_h:= \{ \mathbf{u}\in \mathbf{L}^2( \EO{\Omega}) :  \mathbf{u}|^{}_T\in [\mathbb{P}_0(T)]^{d}\ \forall T \in \mathscr{T}_h\}$. 

The discrete counterpart of \EO{the continuous problem} \eqref{def:weak_ocp}--\eqref{eq:weak_pde} thus reads as follows: Find min $J(\mathbf{y}_h,\mathbf{u}_h)$ subject to the discrete state equations
\begin{equation}\label{eq:discrete_state_eq}
\left\{
\begin{array}{rcll}
a(\mathbf{y}_h,\mathbf{v}_h) + b(\mathbf{v}_h,p_h) & = & (\mathbf{u}_h,\mathbf{v}_h)_{\mathbf{L}^2(\Omega)} & \quad \forall \mathbf{v}_h \in \mathbf{V}_h, \\
b(\mathbf{y}_h,q_h) & = & 0 & \quad \forall q_h \in Q_h,
\end{array}
\right.
\end{equation}
and the discrete control constraints $\mathbf{u}_{h} \in \mathbb{U}_{ad}^{h}$. Standard arguments reveal the existence of a unique optimal pair $(\EO{(\bar{\mathbf{y}}_{h}, \bar{p}_h)},\bar{\mathbf{u}}_{h})$. In addition,  the pair $(\EO{(\bar{\mathbf{y}}_{h}, \bar{p}_h)},\bar{\mathbf{u}}_{h})$ is optimal for the aforementioned discrete optimal control problem if and only if $\EO{(\bar{\mathbf{y}}_{h}, \bar{p}_h)}$ solves \eqref{eq:discrete_state_eq}, with $\mathbf{u}_{h}$ replaced by $\bar{\mathbf{u}}_{h}$, and $\bar{\mathbf{u}}_h$ satisfies the variational inequality
\begin{equation}\label{eq:discrete_variational_ineq}
(\bar{\mathbf{z}}_h + \lambda \bar{\mathbf{u}}_h,\mathbf{u}_h - \bar{\mathbf{u}}_h)_{\mathbf{L}^{2}(\Omega)} \geq 0  \quad \forall \mathbf{u}_h \in \mathbb{U}_{ad}^h,
\end{equation}
where $(\bar{\mathbf{z}}_h, \bar{r}_h)$ solves
\begin{equation}\label{eq:discrete_adj_eq}
\left\{\begin{array}{rcll}
a(\mathbf{w}_h,\mathbf{z}_h) - b(\mathbf{w}_h,r_h)&=&\displaystyle{\sum_{t\in\mathcal{Z}}}\langle ({\mathbf{y}}_h\FF{(t)} - \mathbf{y}_{t})\delta_{t},\mathbf{w}_h \rangle \quad& \forall \mathbf{w}_h \in \mathbf{V}_h, \\
b(\mathbf{z}_h,s_h) &=&0 \quad & \forall  s_h \in Q_h,
\end{array}
\right.\hspace{-0.4cm}
\end{equation}
with ${\mathbf{y}}_h$ replaced by $\bar{\mathbf{y}}_h$. 

\subsubsection{Auxiliary Problems}

We introduce two auxiliary problems that will be instrumental to derive error estimates for the proposed discrete scheme.

The first problem reads as follows: Find $(\hat{\mathbf{y}}_h,\hat{p}_h)\in \mathbf{V}_h\times Q_h$ such that
\begin{align}\label{eq:state_hat}
\begin{cases}
\begin{array}{rcll}
a(\hat{\mathbf{y}}_h,\mathbf{v}_h)+b(\mathbf{v}_h,\hat{p}_h)&=&(\bar{\mathbf{u}},\mathbf{v}_h)_{{\mathbf{L}}^2(\Omega)}\quad& \quad \forall \mathbf{v}_h\in \mathbf{V}_h,\\
b(\hat{\mathbf{y}}_h,q_h) &=&0\quad & \quad \forall  q_h\in Q_h.
\end{array}
\end{cases}
\end{align}
The second auxiliary problem is: Find $(\hat{\mathbf{z}}_h,\hat{r}_h)\in \mathbf{V}_h\times Q_h$ such that
\begin{equation}\label{eq:adjoint_hat}
\left\{\begin{array}{rcll}
a(\mathbf{w}_h,\hat{\mathbf{z}}_h) - b(\mathbf{w}_h,\hat{r}_h)&=&\displaystyle{\sum_{t\in\mathcal{Z}}}\langle (\hat{\mathbf{y}}_h\FF{(t)} - \mathbf{y}_{t})\delta_{t},\mathbf{w}_h \rangle \quad& \quad \forall  \mathbf{w}_h \in \mathbf{V}_h,\\
b(\hat{\mathbf{z}}_h,s_h) &=&0 \quad &\quad \forall  s_h \in Q_h.
\end{array}
\right.\hspace{-0.45cm}
\end{equation}

Before providing error estimates, we present the following auxiliary result.
\begin{lemma}[Discrete pointwise stability]\label{lemma:pointwise_stability_Stokes}
Let $(\boldsymbol\xi,\theta)\in \mathbf{H}_0^1(\Omega)\times L^2(\Omega)/\mathbb{R}$ be the solution to
\begin{align}\label{eq:system_lemma}
\begin{cases}
\begin{array}{rcll}
a(\boldsymbol\xi,\mathbf{v})+b(\mathbf{v},\theta)&=&(\bar{\mathbf{u}}-\bar{\mathbf{u}}_h,\mathbf{v})_{{\mathbf{L}}^2(\Omega)}\quad&\forall \mathbf{v}\in \mathbf{H}_0^1(\Omega),\\
b(\boldsymbol\xi,q) &=&0\quad &\forall q\in L^2(\Omega)/\mathbb{R},
\end{array}
\end{cases}
\end{align}
and let $(\boldsymbol\xi_h,\theta_h)\in \mathbf{V}_h\times Q_h$ be its Galerkin approximation on the basis of the discrete spaces \eqref{def:discrete_spaces_P1B} or \eqref{def:discrete_spaces_TH}. Then,
\begin{equation*}
\|\boldsymbol\xi_h\|_{\mathbf{L}^\infty(\Omega)}
\lesssim
\|\bar{\mathbf{u}}-\bar{\mathbf{u}}_h\|_{\mathbf{L}^2(\Omega)},
\end{equation*}
where the hidden constant is independent of $(\boldsymbol\xi,\theta)$, $(\boldsymbol\xi_h,\theta_h)$, $\bar{\mathbf{u}}$, $\bar{\mathbf{u}}_h$, and $h$.
\end{lemma}
{\it Proof.}
We begin the proof by noticing that, since $\bar{\mathbf{u}} - \bar{\mathbf{u}}_h \in \mathbf{L}^2(\Omega)$, then $(\boldsymbol\xi,\theta)\in \mathbf{H}^2(\Omega) \cap \mathbf{H}_0^1(\Omega)\times H^1(\Omega) \cap L^2(\Omega)/\mathbb{R}$. This results follows from Proposition \ref{pro:H2regularity}.
Let us denote by $\mathbf{I}_h:\mathbf{C}(\bar{\Omega})\rightarrow \mathbf{V}_h$ the Lagrange interpolation operator. An application of the triangle inequality in conjunction with a standard inverse estimate yield
\begin{multline}
\|\boldsymbol\xi_h\|_{\mathbf{L}^\infty(\Omega)}
\lesssim 
\|\boldsymbol\xi\|_{\mathbf{L}^\infty(\Omega)}+\|\boldsymbol\xi-\mathbf{I}_h\boldsymbol\xi\|_{\mathbf{L}^\infty(\Omega)}+h^{-d/2}\|\mathbf{I}_h\boldsymbol\xi-\boldsymbol\xi_h\|_{\mathbf{L}^2(\Omega)}\\ 
\lesssim
\|\boldsymbol\xi\|_{\mathbf{L}^\infty(\Omega)}+\|\boldsymbol\xi-\mathbf{I}_h\boldsymbol\xi\|_{\mathbf{L}^\infty(\Omega)}\\
+h^{-d/2}\left(\|\mathbf{I}_h\boldsymbol\xi-\boldsymbol\xi\|_{\mathbf{L}^2(\Omega)}+\|\boldsymbol\xi-\boldsymbol\xi_h\|_{\mathbf{L}^2(\Omega)}\right).
\label{eq:ineq_xi}
\end{multline}
To control the first term on the right hand side of \eqref{eq:ineq_xi} we invoke the continuous Sobolev embedding $\mathbf{H}^{2}\hookrightarrow\mathbf{C}(\bar{\Omega})$ and the regularity estimate \eqref{eq:H2regularity} to arrive at
\begin{equation*}\label{eq:xi_embedding}
\|\boldsymbol\xi\|_{\mathbf{L}^\infty(\Omega)}
\lesssim
\|\boldsymbol\xi\|_{\mathbf{H}^{2}(\Omega)}
\lesssim
\|\bar{\mathbf{u}}-\bar{\mathbf{u}}_h\|_{\mathbf{L}^2(\Omega)}.
\end{equation*}
For the remaining terms in the right hand side of \eqref{eq:ineq_xi} we utilize standard interpolation error estimates for $\mathbf{I}_{h}$ and error estimates for the finite element approximation of problem \eqref{eq:system_lemma}. This concludes the proof.
\qed

\subsubsection{Error Estimates}
To perform an a priori error analysis, it is useful to introduce the $\mathbf{L}^2(\Omega)$--orthogonal projection onto $ [\mathbb{P}_0(\mathscr{T}_h)]^d$, which is defined by
%
\[
\Pi_{\mathbf{L}^2}: \mathbf{L}^2(\Omega) \rightarrow [\mathbb{P}_0(\mathscr{T}_h)]^d,
\quad
 \left(\Pi_{\mathbf{L}^2}\mathbf{v}\right)|^{}_T:=\frac{1}{|T|}\left(\int_T v_1(x) \: \mathrm{d}x,\dots,\int_T v_d(x) \: \mathrm{d}x\right),
\]
where $T \in \mathscr{T}_h$. Note that, in view of the weighted Poincar\'e inequality of Theorem \ref{thm:weighted_poincare}, $\Pi_{\mathbf{L}^2}$ is well defined over $\mathbf{H}_0^1(\rho,\Omega)$.

\EO{To obtain error estimates for $\Pi_{\mathbf{L}^2}$, we proceed with the help of regularity estimates for a suitable Green's function. For $x \in \Omega$, we define the Green's function $\mathcal{G}: \Omega \times \Omega \rightarrow \mathbb{R}$ as the solution (in the sense of distributions) to 
\[
-\Delta_y \mathcal{G} = \delta (x-y), 
\quad y \in \Omega, 
\quad
\mathcal{G}(x,y) = 0, 
\quad y \in \partial \Omega.
\]
The following properties of the Green's function are essential in what follows.

\begin{proposition}[Properties of $\mathcal{G}$]
Let $\Omega \subset \mathbb{R}^d$ be a Lipschitz polytope. The Green's function $\mathcal{G}$ satisfies the following properties:
\begin{enumerate}
\item $\nabla \mathcal{G}(x,\cdot) \in L^{\frac{d}{d-1},\infty}(\Omega)$.
\item $\mathcal{G}(x,\cdot) \in W^{1,\frac{d}{d-1}}(\Omega \setminus B_R)$ and $\mathcal{G}(x,\cdot) \in W^{2,1}(\Omega \setminus B_R)$, where $B_R$ denotes the ball of radius $R$ and center $x \in \Omega$. In addition, we have
\begin{equation}
 \| \nabla \mathcal{G}(x,\cdot) \|_{L^{\frac{d}{d-1}}(\Omega \setminus B_R)} 
 +
 \| D^2 \mathcal{G}(x,\cdot) \|_{L^1(\Omega \setminus B_R)} \lesssim | \log R^{-1}|.
 \label{eq:Green_W21}
\end{equation}
\end{enumerate}
\label{pro:Green}
\end{proposition}
{\it Proof.} A proof of $\nabla \mathcal{G} \in L^{d/(d-1),\infty}(\Omega)$ can be found in \cite[Theorem 1.1]{MR657523} and \cite[Theorem 4.1]{MR2341783}, for $d=3$, and \cite[Theorem 1.1]{MR2763343} and \cite[Section 6]{MR1354111}, for $d=2$. Here, $L^{s,\infty}(\Omega)$, with $0 < s < \infty $, corresponds to the so--called weak Lebesgue space; see definition 1.15 in \cite{MR2445437}. A proof of the estimate \eqref{eq:Green_W21} can be found in \cite[Theorem 1]{kopteva}.
\qed

We now provide an error estimate for the $\Pi_{L^2}$ projection.

\begin{theorem}
Let $\Omega$ be \FF{an open,} bounded and Lipschitz polytope. Let $x_0 \in \Omega$ and let $u$ be the weak solution to the Poisson problem $-\Delta u = \delta_{x_0}$ in $\Omega$ and $u = 0$ on $\partial \Omega$. Thus, we have the error estimate
\begin{equation}
\label{eq:error_estimate_projection}
\| u - \Pi_{L^2} u \|_{L^2(\Omega)} \lesssim |\log h| h^{2-d/2},
\end{equation}
with a hidden constant that is independent of $h$.
\label{thm:error_estimate_projection}
\end{theorem}
{\it Proof.} Let us consider a partition of $\mathscr{T}_h$ into the sets
$
\mathcal{N}_{x_0} := \{ T \in \mathscr{T}_h: x_0 \in \mathcal{N}_T \}
$
and $\mathcal{M}_{x_0}:= \mathscr{T}_h \setminus \mathcal{N}_{x_0}$ and utilize the regularity results of Proposition \ref{pro:Green} accordingly; we recall that $\mathcal{N}_T$ is defined in \eqref{def:patch}. We proceed in two steps.

\textit{Step 1}. Let $T \in \mathcal{M}_{x_0}$. Standard approximation results yields the estimate
\[
\| u - \Pi_{L^2} u \|_{L^2(T)} \lesssim h_T^{2-d/2} \| \nabla u \|_{L^{\frac{d}{d-1}}(\mathcal{N}_T)}.
\]
Since the mesh $\mathscr{T}_h$ is quasi--uniform, we thus obtain
\begin{equation}
\begin{aligned}
\sum_{T \in \mathcal{M}_{x_0}} \| u - \Pi_{L^2} u \|^2_{L^2(T)} & \lesssim h^{2-d/2}  \sum_{T \in \mathcal{M}_{x_0}}  \| \nabla u \|_{L^{\frac{d}{d-1}}(\mathcal{N}_T)} 
\\
& \lesssim |\log h|h^{2-d/2},
\end{aligned}
\label{eq:proyection_estimate_1}
\end{equation}
where, to obtain the last estimate, we have used the finite overlapping property of stars and \eqref{eq:Green_W21} to conclude that $\| \nabla  u \|_{L^{d/(d-1)}(\Omega \setminus B_R)} \lesssim | \log h|$, with $B_R$, in this case, being the largest ball such that $B_R \subset \cup \{T \in \mathscr{T}_h: x_0 \in T\}$.
 
\textit{Step 2.} Let $T \in \mathcal{N}_{x_0}$. Observe that
\begin{align*}
\| u - \Pi_{L^2} u \|^2_{L^2(T)}  
& =  \int_{T}\left|\frac{1}{|T|} \int_{T} \left( u(x) - u(y) \right) \mathrm{d} y\right|^2 \mathrm{d}x
\\
& =   \int_{T}\left|\frac{1}{|T|} \int_{T} \int_0^1 \nabla u(tx+(1-t)y) \cdot (x-y) \mathrm{d}t \mathrm{d} y\right|^2\mathrm{d}x
\\
& \lesssim h_T^2 \int_{T} \left|\frac{1}{|T|} \int_{T} \nabla u(z) \mathrm{d} z \right|^2 \mathrm{d}x
 \lesssim  
h_T^{2-d} \left|\int_{T} \nabla u(z) \mathrm{d} z \right|^2.
\end{align*}

Let us use now the fact that $\nabla u \in L^{d/(d-1),\infty}(\Omega)$. Since $\nabla u \in L^{s,\infty}(\Omega)$ for $s \leq d/(d-1)$, we have that the distribution function of $|\nabla u|$, i.e., $d_{\nabla u}(\gamma):= | \{ x \in \Omega : |\nabla u (x)| > \gamma \} |$, satisfies $d_{\nabla u} (\gamma) \lesssim \gamma^{-s}$ for every $\gamma > 0$; see \cite[(1.1.8) and (1.1.9)]{MR2445437}. Here, $| A |$ denotes the Lebesgue measure of a Lebesgue measurable subset $A$ of $\mathbb{R}^d$. We recall that
\[
\| \nabla u \|_{L^{s,\infty}(\Omega)}:= \inf \{ C >0: d_{\nabla u}(\gamma) \leq C^s \gamma^{-s} \}.
\]

Let $C_T := \| \nabla u\|_{L^{s,\infty}(T)}$. Invoke \cite[Proposition 1.1.4]{MR2445437} to obtain, for $s \in (1, d/(d-1)]$, that
\begin{align*}
\int_T |\nabla u| \mathrm{d}x 
= \int_{0}^{\infty} d_{\nabla u}(\gamma) \mathrm{d} \gamma
 \leq \int_0^{\sigma} |T| \mathrm{d} \gamma + \int_{\sigma}^{\infty }  C_T^{s} \gamma^{-s} \mathrm{d} \gamma  
=  \sigma |T| -  C_T^s \frac{ \sigma^{1-s}}{1-s}.
\end{align*}
Set $\sigma = C_T [|T| (s-1)]^{-1/s}$ to arrive at $\int_T |\nabla u| \mathrm{d}x  \lesssim C_T h_T^{d(1-1/s)}$. Set now $s = d/(d-1)$. This yields $\int_T |\nabla u| \mathrm{d}x  \lesssim C_T h_T$. Consequently, 
\begin{equation*}
\| u - \Pi_{L^2} u \|^2_{L^2(T)} 
\lesssim 
h_T^{2-d} \left|\int_{T} \nabla u(z) \mathrm{d} z \right|^2
C_T^2
\lesssim h^{4-d} C_T^2,
\end{equation*}
with a hidden constant independent of $u$ and $h$. 
The previous estimate implies 
\begin{equation}
\sum_{T \in \mathcal{N}_{x_0}} \| u - \Pi_{L^2} u \|^2_{L^2(\Omega)} \lesssim h^{4-d} \sum_{T \in \mathcal{N}_{x_0}} C_T^2 \lesssim \max_{T \in \mathcal{N}_{x_0}} \{ C_T^2 \} h^{4-d},
\label{eq:proyection_estimate_2}
\end{equation}
where we have used that, since $\{ \mathscr{T}_h \}$ is shape--regular, $\# \mathcal{N}_{x_0}$ is uniformly bounded. This also yields the bound $\max_{T \in \mathcal{N}_{x_0}} \{ C_T^2 \} \lesssim 1$, with a hidden constant independent of $h$.

A collection of the bounds \eqref{eq:proyection_estimate_1} and \eqref{eq:proyection_estimate_2} yields the desired estimate \eqref{eq:error_estimate_projection}. This concludes the proof. 
\qed

\begin{remark}[Green's function for Stokes problem and error estimate for $\Pi_{\mathbf{L}^2}$]
The Dirichlet Green's matrix associated with the Stokes system can be defined as follows: $\mathcal{G} = (\mathcal{Q},\mathcal{R})$, where $\mathcal{Q}$ is a $d \times d$ matrix--valued function and $\mathcal{R}$ is a $1 \times d$ vector--valued function that satisfy system (7.8) in \cite{MR2763343}. Theorem 7.1 in \cite{MR2763343} yields, for $d \in \{ 2,3 \}$, 
\begin{equation}
\label{eq:Q_weak_Lp}
\nabla \mathcal{Q}(x,\cdot) \in \mathbf{L}^{\frac{d}{d-1},\infty}(\Omega);
\end{equation}
see also \cite[Theorem 1.1]{MR3320459} and \cite[Theorem 2.4]{MR3877495}. Estimate \eqref{eq:Green_W21} does not appear directly in the literature to our knowledge. In the case of Poisson's problem, we stress that \eqref{eq:Green_W21} can be found in \cite[Lemma 2]{kopteva}. We shall thus assume estimate \eqref{eq:Green_W21} in the case of the Stokes system. With this assumption and the regularity estimate \eqref{eq:Q_weak_Lp} at hand, the arguments elaborated in Theorem \ref{thm:error_estimate_projection} yield the error estimate
\begin{equation}
\label{eq:error_estimate_projection_z}
\| \bar{\mathbf{z}} - \Pi_{\mathbf{L}^2} \bar{\mathbf{z}}  \|_{\mathbf{L}^2(\Omega)} 
\lesssim
|\log h| h^{2-d/2},
\end{equation}
where the hidden constant depends on the $\mathbf{L}^{d/(d-1),\infty}(\Omega)$-norm of $\bar{\mathbf{z}}$ but, more importantly, is independent of $h$. In view of the projection formula \eqref{eq:projection_formula}, \cite[Theorem A.1]{MR567696}, and the definition of the weak Lebesgue space $\mathbf{L}^{s,\infty}(\Omega)$, with $0<s<\infty$, we can conclude that the regularity estimates needed to obtain \eqref{eq:error_estimate_projection_z} are inherited to the control variable $\bar{\mathbf{u}}$. Consequently,
\begin{equation}
\label{eq:error_estimate_projection_u}
\| \bar{\mathbf{u}} - \Pi_{\mathbf{L}^2} \bar{\mathbf{u}}  \|_{\mathbf{L}^2(\Omega)} 
\lesssim
|\log h| h^{2-d/2},
\end{equation}
with, again, a hidden constant independent of $h$.
\end{remark}

We are ready to prove the main result of this section.}

\begin{theorem}[Rates of convergence for $\bar{\mathbf{u}}$]
\label{thm:control_rates}
Let $(\bar{\mathbf{y}},\bar{p},\bar{\mathbf{z}},\bar{r},\bar{\mathbf{u}})\in \mathbf{H}_0^1(\Omega)\times L^2(\Omega)/\mathbb{R}\\ \times \mathbf{H}_0^1(\rho,\Omega) \times L^2(\rho,\Omega)/\mathbb{R}\times \mathbb{U}_{ad}$ be the solution to the optimality system \eqref{eq:weak_pde}, \eqref{eq:adj_eq}, and \eqref{eq:variational_ineq} and $(\bar{\mathbf{y}}_h,\bar{p}_h,\bar{\mathbf{z}}_h,\bar{r}_h,\bar{\mathbf{u}}_h)\in \mathbf{V}_h\times Q_h\times\mathbf{V}_h\times Q_h\times \mathbb{U}_{ad}^h$ its numerical approximation given as the solution to \eqref{eq:discrete_state_eq}--\eqref{eq:discrete_adj_eq}. \EO{Then,
\begin{multline}\label{eq:global_reliability_2}
\| \bar{\mathbf{u}} - \bar{\mathbf{u}}_h  \|_{\mathbf{L}^2(\Omega)}
\lesssim
h^{2-d/2}\sum_{t\in\mathcal{Z}} |\bar{\mathbf{y}}(t) - \mathbf{y}_{t}| \|\delta_{t}\|_{\mathcal{M}(\Omega)} \\
+ h^{2-d/2}|\log h| 
+ h|\log h|^{3}\| \bar{\mathbf{u}} \|_{\mathbf{L}^{\vartheta}(\Omega)},
\end{multline}
with $\vartheta > d$. The} hidden constant is independent of the continuous and discrete solutions, the size of the elements in the mesh $\mathscr{T}_h$, and $\#\mathscr{T}_h$. The constant, however, blows up as \FF{$\lambda\downarrow 0$}. 
\end{theorem}
{\it Proof.}
We proceed in four steps.

\underline{Step 1.} Let us consider $\mathbf{u}=\bar{\mathbf{u}}_h$ in \eqref{eq:variational_ineq} and $\mathbf{u}_h=\Pi_{\mathbf{L}^2}\bar{\mathbf{u}}$ in \eqref{eq:discrete_variational_ineq}. Adding the obtained inequalities we arrive at
\begin{equation}\label{eq:ineq_1}
\lambda\|\bar{\mathbf{u}}-\bar{\mathbf{u}}_h\|_{\mathbf{L}^2(\Omega)}^2
\leq
(\bar{\mathbf{z}}-\bar{\mathbf{z}}_h,\bar{\mathbf{u}}_h-\bar{\mathbf{u}})_{\mathbf{L}^2(\Omega)}+(\bar{\mathbf{z}}_h+\lambda\bar{\mathbf{u}}_h,\Pi_{\mathbf{L}^2}\bar{\mathbf{u}}-\bar{\mathbf{u}})_{\mathbf{L}^2(\Omega)}.
\end{equation}

\underline{Step 2.} The goal of this step is to bound the term $(\bar{\mathbf{z}}-\bar{\mathbf{z}}_h,\bar{\mathbf{u}}_h-\bar{\mathbf{u}})_{\mathbf{L}^2(\Omega)}$. To accomplish this task, we add and subtract the auxiliary term $\hat{\mathbf{z}}_h$, where $(\hat{\mathbf{z}}_h,\hat r_h)$ corresponds to the solution to \eqref{eq:adjoint_hat}, to obtain
\begin{align*}
(\bar{\mathbf{z}}_h-\bar{\mathbf{z}},\bar{\mathbf{u}}-\bar{\mathbf{u}}_h)_{\mathbf{L}^2(\Omega)}
&=
(\bar{\mathbf{z}}_h-\hat{\mathbf{z}}_h,\bar{\mathbf{u}}-\bar{\mathbf{u}}_h)_{\mathbf{L}^2(\Omega)}+(\hat{\mathbf{z}}_h-\bar{\mathbf{z}},\bar{\mathbf{u}}-\bar{\mathbf{u}}_h)_{\mathbf{L}^2(\Omega)}\\
&=:\mathbf{I}+\mathbf{II}.
\end{align*}
Let us concentrate on $\mathbf{I}$. Note that $(\bar{\mathbf{y}}_h-\hat{\mathbf{y}}_h,\bar{p}_h-\hat{p}_h)\in \mathbf{V}_h\times Q_h$ solves
\begin{align}
\begin{cases}\label{eq:state_bar_hat}
\begin{array}{rcll}
a(\bar{\mathbf{y}}_h-\hat{\mathbf{y}}_h,\mathbf{v}_h)+b(\mathbf{v}_h,\bar{p}_h-\hat{p}_h)&=&(\bar{\mathbf{u}}_h-\bar{\mathbf{u}},\mathbf{v}_h)_{{\mathbf{L}}^2(\Omega)}, 
\\
b(\bar{\mathbf{y}}_h-\hat{\mathbf{y}}_h,q_h) &=&0
\end{array}
\end{cases}
\end{align}
for all $\mathbf{v}_h\in \mathbf{V}_h$ and $q_h\in Q_h$, and that $(\bar{\mathbf{z}}_h-\hat{\mathbf{z}}_h,\bar{r}_h-\hat{r}_h)\in \mathbf{V}_h\times Q_h$ solves 
\begin{align}
\label{eq:adjoint_bar_hat}
\begin{cases}
\begin{array}{rcll}
a(\mathbf{w}_h,\bar{\mathbf{z}}_h-\hat{\mathbf{z}}_h) - b(\mathbf{w}_h,\bar{r}_h-\hat{r}_h)&=&\displaystyle{\sum_{t\in\mathcal{Z}}}\langle (\bar{\mathbf{y}}_h\FF{(t)} - \hat{\mathbf{y}}_h\FF{(t)})\delta_{t},\mathbf{w}_h\rangle,
\\
b(\bar{\mathbf{z}}_h-\hat{\mathbf{z}}_h,s_h) &=&0
\end{array}
\end{cases}
\end{align}
for all $\mathbf{w}_h \in \mathbf{V}_h$ and $s_h \in Q_h$. Set $\mathbf{w}_h=\bar{\mathbf{y}}_h-\hat{\mathbf{y}}_h\in \mathbf{V}_h$ in \eqref{eq:adjoint_bar_hat} and $\mathbf{v}_h=\bar{\mathbf{z}}_h-\hat{\mathbf{z}}_h\in \mathbf{V}_h$ in \eqref{eq:state_bar_hat} to conclude 
\begin{equation}
\label{eq:I}
\mathbf{I}=-\sum_{t\in\mathcal{Z}}|\bar{\mathbf{y}}_h(t)-\hat{\mathbf{y}}_h(t)|^2\leq 0.
\end{equation}
To estimate $\mathbf{II}$ we proceed as follows. First, we use Young's inequality to obtain
\begin{equation*}
\mathbf{II}
\leq
\frac{\lambda}{4}\|\bar{\mathbf{u}}-\bar{\mathbf{u}}_h\|_{\mathbf{L}^2(\Omega)}^2+\frac{1}{\lambda}\|\hat{\mathbf{z}}_h-\bar{\mathbf{z}}\|_{\mathbf{L}^2(\Omega)}^2.
\end{equation*}
Second, to estimate the term $\|\hat{\mathbf{z}}_h-\bar{\mathbf{z}}\|_{\mathbf{L}^2(\Omega)}^2$, we introduce $(\tilde{\mathbf{z}}_h,\tilde{r}_h)\in \mathbf{V}_h\times Q_h$ as the solution to
\begin{equation}\label{eq:adjoint_tilde}
\left\{\begin{array}{rcll}
a(\mathbf{w}_h,\tilde{\mathbf{z}}_h) - b(\mathbf{w}_h,\tilde{r}_h)&=&\displaystyle{\sum_{t\in\mathcal{Z}}}\langle (\bar{\mathbf{y}}\FF{(t)} - \mathbf{y}_{t})\delta_{t},\mathbf{w}_h \rangle \quad&\forall \mathbf{w}_h \in \mathbf{V}_h,
\\
b(\tilde{\mathbf{z}}_h,s_h) &=&0 \quad &\forall s_h \in Q_h.
\end{array}
\right.\hspace{-0.45cm}
\end{equation}
Third, add and subtract $\tilde{\mathbf{z}}_h$ and use the triangle inequality to arrive at
\begin{equation}\label{eq:ineq_3.5}
\|\hat{\mathbf{z}}_h-\bar{\mathbf{z}}\|_{\mathbf{L}^2(\Omega)}^2
\leq
2\|\hat{\mathbf{z}}_h-\tilde{\mathbf{z}}_{h}\|_{\mathbf{L}^2(\Omega)}^2+2\|\tilde{\mathbf{z}}_h-\bar{\mathbf{z}}\|_{\mathbf{L}^2(\Omega)}^2.
\end{equation}

We now analyze $\|\tilde{\mathbf{z}}_h-\bar{\mathbf{z}}\|_{\mathbf{L}^2(\Omega)}$. Since $(\tilde{\mathbf{z}}_h,\tilde{r}_h)$ corresponds to the Galerkin approximation of $(\bar{\mathbf{z}},\bar{r})$, an application of Lemma \ref{lemma:first_estimate} yields
\begin{equation}\label{eq:ineq_4}
\|\tilde{\mathbf{z}}_h-\bar{\mathbf{z}}\|_{\mathbf{L}^2(\Omega)}
\lesssim 
\FF{h^{2-d/2}\sum_{t\in\mathcal{Z}} |\bar{\mathbf{y}}(t) - \mathbf{y}_{t}|\|\delta_{t}\|_{\mathcal{M}(\Omega)}}.
\end{equation}
It thus remains to estimate $\|\hat{\mathbf{z}}_h-\tilde{\mathbf{z}}_h\|_{\mathbf{L}^2(\Omega)}$. To accomplish this, we invoke the weighted Poincar\'e inequality of Theorem \ref{thm:weighted_poincare} to obtain the estimate
\begin{equation*}
\|\hat{\mathbf{z}}_h-\tilde{\mathbf{z}}_h\|_{\mathbf{L}^2(\Omega)}^2
\lesssim
\|\nabla(\hat{\mathbf{z}}_h-\tilde{\mathbf{z}}_h)\|_{\mathbf{L}^2(\rho,\Omega)}^2.
\end{equation*}
This, in view of the stability of the discrete Stokes system in weighted spaces \cite[Theorem 4.1]{2019arXiv190500476D}, allows us to obtain 
\begin{equation}\label{eq:ineq_6}
\|\hat{\mathbf{z}}_h-\tilde{\mathbf{z}}_h\|_{\mathbf{L}^2(\Omega)}^2
\lesssim
\|\nabla(\hat{\mathbf{z}}_h-\tilde{\mathbf{z}}_h)\|_{\mathbf{L}^2(\rho,\Omega)}^2 
\lesssim
\| \hat{\mathbf{y}}_h - \bar{\mathbf{y}} \|_{\mathbf{L}^{\infty}(\Omega)}^{2}.
\end{equation}
Now, let us recall that $\hat{\mathbf{y}}_{h}$ is the Galerkin approximation of $\bar{\mathbf{y}}$. In addition, since $d\in\{2,3\}$, $\Omega$ is a convex polytope, and $\bar{\mathbf{u}}\in\mathbf{L}^\infty(\Omega)$, we have $\bar{\mathbf{y}}\in \FF{\mathbf{W}}^{1,\infty}(\Omega)$ \cite{MR2228352} (see also  \cite{MR2945141,MR3422453}). Therefore, the pointwise error estimates of \cite[Theorem 4.1]{MR995211}, combined with the weighted estimates of \cite{MR1880723} for $d=3$, yield
\begin{equation}
\label{eq:ineq_7}
\begin{aligned}
\| \bar{\mathbf{y}} - \hat{\mathbf{y}}_h \|_{\mathbf{L}^{\infty}(\Omega)}^{2}
& \lesssim
h^2|\log h|^{6}\left (\|\nabla\bar{\mathbf{y}}\|_{\mathbf{L}^{\infty}(\Omega)}^{2} + \|\bar{p}\|_{L^{\infty}(\Omega)}^{2} \right) 
\\
& \lesssim \EO{h^2|\log h|^{6} \| \bar{\mathbf{u}} \|_{\mathbf{L}^{\vartheta}(\Omega)}},
\end{aligned}
\end{equation}
\EO{with $\vartheta>d$. To obtain the last estimate, we have used \cite[estimate (1.17)]{MR3422453}.} Recall that $\mathbf{I} \leq 0$. We thus replace \eqref{eq:ineq_7} into \eqref{eq:ineq_6} and combine the obtained estimate with \eqref{eq:ineq_4} to conclude that 
\begin{multline}\label{eq:ineq_8}
(\bar{\mathbf{z}}-\bar{\mathbf{z}}_h,\bar{\mathbf{u}}_h-\bar{\mathbf{u}})_{\mathbf{L}^2(\Omega)} 
\leq
\frac{\lambda}{4}\|\bar{\mathbf{u}}-\bar{\mathbf{u}}_h\|_{\mathbf{L}^2(\Omega)}^2 \\
\EO{+ \frac{C_1}{\lambda}\left(h^{4-d}\sum_{t\in\mathcal{Z}} |\bar{\mathbf{y}}(t) - \mathbf{y}_{t}|^2\|\delta_{t}\|_{\mathcal{M}(\Omega)}^2 + C_2 h^2|\log h|^{6}\| \bar{\mathbf{u}} \|_{\mathbf{L}^{\vartheta}(\Omega)}^2\right),}
\end{multline}
\FF{with $\vartheta>d$, and} $C_1$ and $C_2$ being positive constants.

\underline{Step 3.} The goal of this step is to bound the remaining term in \eqref{eq:ineq_1}. Note that, by adding and subtracting the term $\lambda \bar{\mathbf{u}}$ and the adjoint variables $\bar{\mathbf{z}}$ and $\hat{\mathbf{z}}_h$, we obtain the following identity: 
\begin{multline*}
(\bar{\mathbf{z}}_h+\lambda\bar{\mathbf{u}}_h,\Pi_{\mathbf{L}^2}\bar{\mathbf{u}}-\bar{\mathbf{u}})_{\mathbf{L}^2(\Omega)}
=
(\bar{\mathbf{z}} + \lambda\bar{\mathbf{u}},\Pi_{\mathbf{L}^2}\bar{\mathbf{u}}-\bar{\mathbf{u}})_{\mathbf{L}^2(\Omega)} + (\bar{\mathbf{z}}_h - \hat{\mathbf{z}}_h ,\Pi_{\mathbf{L}^2}\bar{\mathbf{u}}-\bar{\mathbf{u}})_{\mathbf{L}^2(\Omega)} \\+ (\hat{\mathbf{z}}_h - \bar{\mathbf{z}} ,\Pi_{\mathbf{L}^2}\bar{\mathbf{u}}-\bar{\mathbf{u}})_{\mathbf{L}^2(\Omega)}
 + \lambda(\bar{\mathbf{u}}_h - \bar{\mathbf{u}},\Pi_{\mathbf{L}^2}\bar{\mathbf{u}}-\bar{\mathbf{u}})_{\mathbf{L}^2(\Omega)} =: \mathbf{III}_1+\mathbf{III}_2+ \mathbf{III}_3+ \mathbf{III}_4.
\end{multline*}

We now bound the terms $\mathbf{III}_1$, $\mathbf{III}_2$, $\mathbf{III}_3$, and $\mathbf{III}_4$.
To estimate $\mathbf{III}_{1}$ we \EO{invoke the error estimates \eqref{eq:error_estimate_projection_z} and \eqref{eq:error_estimate_projection_u}. We thus proceed} as follows:
\begin{multline}
\label{eq:term_I}
\mathbf{III}_1
= 
(\bar{\mathbf{z}} + \lambda\bar{\mathbf{u}}-\Pi_{\mathbf{L}^2}\bar{\mathbf{z}} - \lambda\Pi_{\mathbf{L}^2}\bar{\mathbf{u}},\Pi_{\mathbf{L}^2}\bar{\mathbf{u}}-\bar{\mathbf{u}})_{\mathbf{L}^2(\Omega)} \\ 
\leq 
\frac{1}{2} \|\Pi_{\mathbf{L}^2}\bar{\mathbf{u}}-\bar{\mathbf{u}}\|_{\mathbf{L}^2(\Omega)}^2+\frac{1}{2}\|\bar{\mathbf{z}}-\Pi_{\mathbf{L}^2}\bar{\mathbf{z}}\|_{\mathbf{L}^2(\Omega)}^2
\lesssim
|\log h|^2 h^{4-d}.
\end{multline}

Now we estimate $\mathbf{III}_2$. \EO{We begin with a simple application of the Cauchy--Schwarz inequality and the Poincar\'e inequality of Theorem \ref{thm:weighted_poincare} to write
\begin{equation*}
\begin{aligned}
\mathbf{III}_2
& \leq 
\|\Pi_{\mathbf{L}^2}\bar{\mathbf{u}}-\bar{\mathbf{u}}\|_{\mathbf{L}^2(\Omega)}
\|\bar{\mathbf{z}}_h-\hat{\mathbf{z}}_h \|_{\mathbf{L}^2(\Omega)}
\\
& \lesssim \|\Pi_{\mathbf{L}^2}\bar{\mathbf{u}}-\bar{\mathbf{u}}\|_{\mathbf{L}^2(\Omega)} \|\nabla(\bar{\mathbf{z}}_h-\hat{\mathbf{z}}_h)\|_{\mathbf{L}^2(\rho,\Omega)}.
\end{aligned}
\end{equation*}
Invoking} the stability of the discrete Stokes system in weighted spaces of \cite[Theorem 4.1]{2019arXiv190500476D}, we obtain
\begin{equation}\label{eq:stability_bar_hat}
\|\nabla(\bar{\mathbf{z}}_h-\hat{\mathbf{z}}_h)\|_{\mathbf{L}^2(\rho,\Omega)}
\lesssim
\|\bar{\mathbf{y}}_h-\hat{\mathbf{y}}_h\|_{\mathbf{L}^\infty(\Omega)}.
\end{equation}
To estimate the right hand side of the previous expression we introduce the auxiliary variables $(\hat{\mathbf{y}},\hat{p})\in \mathbf{H}_0^1(\Omega)\times L^2(\Omega)/\mathbb{R}$ as the solution to
\begin{align}\label{eq:state_hat_continuous}
\begin{cases}
\begin{array}{rcll}
a(\hat{\mathbf{y}},\mathbf{v})+b(\mathbf{v},\hat{p})&=&(\bar{\mathbf{u}}_h,\mathbf{v})_{{\mathbf{L}}^2(\Omega)}\quad&\forall  \mathbf{v}\in \mathbf{H}_0^1(\Omega),\\
b(\hat{\mathbf{y}},q) &=&0\quad &\forall  q\in L^2(\Omega)/\mathbb{R}.
\end{array}
\end{cases}
\end{align}
Since the pair $(\bar{\mathbf{y}}-\hat{\mathbf{y}},\bar{p}-\hat{p})$ solves the Stokes system with the term $\bar{\mathbf{u}}-\bar{\mathbf{u}}_h$ in the right hand side of the \emph{momentum equation} and \EO{$(\hat{\mathbf{y}}_h-\bar{\mathbf{y}}_h,\hat{p}_h-\bar{p}_h)$} corresponds to the Galerkin approximation of $(\bar{\mathbf{y}}-\hat{\mathbf{y}},\bar{p}-\hat{p})$, it follows from Lemma \ref{lemma:pointwise_stability_Stokes} that
\begin{equation}\label{eq:ineq_xitau}
\|\bar{\mathbf{y}}_h-\hat{\mathbf{y}}_h\|_{\mathbf{L}^\infty(\Omega)} 
\lesssim
\|\bar{\mathbf{u}}-\bar{\mathbf{u}}_h\|_{\mathbf{L}^2(\Omega)},
\end{equation}
where we have considered, \EO{in the notation of Lemma \ref{lemma:pointwise_stability_Stokes}}, $(\boldsymbol\xi,\theta)=(\bar{\mathbf{y}}-\hat{\mathbf{y}},\bar{p}-\hat{p})\in \mathbf{H}_0^1(\Omega)\times L^2(\Omega)/\mathbb{R}$ and $(\boldsymbol\xi_h,\theta_h)=(\hat{\mathbf{y}}_h-\bar{\mathbf{y}}_h,\hat{p}_h-\bar{p}_h)\in \mathbf{V}_h\times Q_h$. We thus invoke \eqref{eq:stability_bar_hat}, \eqref{eq:ineq_xitau}, \EO{the error estimate \eqref{eq:error_estimate_projection_u}, and Young's inequality to} conclude that
\begin{equation}\label{eq:term_II}
\mathbf{III}_2 
\leq
\EO{\frac{C}{\lambda} h^{4-d} |\log h|^2 + \frac{\lambda}{4} \|\bar{\mathbf{u}} - \bar{\mathbf{u}}_{h}\|_{\mathbf{L}^{2}(\Omega)}^{2}},
\end{equation}
\EO{where $C$ denotes a positive constant that is independent of $h$ and $\lambda$.}

We bound the term $\mathbf{III}_3$ by using Young's inequality and the estimates provided in \eqref{eq:ineq_3.5}--\eqref{eq:ineq_7}. These arguments reveal that
\begin{equation}\label{eq:term_III}
\mathbf{III}_3
\lesssim
\EO{h^{4-d}\sum_{t\in\mathcal{Z}} |\bar{\mathbf{y}}(t) - \mathbf{y}_{t}|^2\|\delta_{t}\|^2_{\mathcal{M}(\Omega)}} 
\EO{+ h^{4-d} |\log h|^2}
 + h^2|\log h|^{6}\FF{\| \bar{\mathbf{u}} \|_{\mathbf{L}^{\vartheta}(\Omega)}^2},
\end{equation}
\FF{with $\vartheta>d$.}
To estimate $\mathbf{III}_4$ we use Young's inequality to \EO{immediately} arrive at
\begin{equation}\label{eq:term_IV}
\mathbf{III}_4 \leq \dfrac{\lambda}{4} \| \bar{\mathbf{u}} - \bar{\mathbf{u}}_{h} \|_{\mathbf{L}^{2}(\Omega)}^{2}  + C \lambda \EO{h^{4-d} |\log h|^2},
\end{equation}
where $C$ denotes a positive constant that is independent of \EO{$h$ and} $\lambda$.

\underline{Step 4.} The proof concludes by gathering \eqref{eq:ineq_1}, \eqref{eq:ineq_8}, \eqref{eq:term_I}, \eqref{eq:term_II}, \eqref{eq:term_III}, and \eqref{eq:term_IV}.
\qed

\begin{remark}[Rates of convergence for $\bar{\mathbf{u}}$]\label{rem:better_control_rates}
\EO{The error estimate of Theorem \ref{thm:control_rates} reads as follows:
\[
\| \bar{\mathbf{u}} - \bar{\mathbf{u}}_h  \|_{\mathbf{L}^2(\Omega)} \lesssim h|\log h|^3,
\quad
\| \bar{\mathbf{u}} - \bar{\mathbf{u}}_h  \|_{\mathbf{L}^2(\Omega)} \lesssim h^{1/2} |\log h|,
\]
for $d=2$ and $d=3$, respectively. The two--dimensional error estimate is nearly--optimal in terms of approximation (nearly because of the presence of the $\log$-term). This rate of convergence is dictated by the regularity properties of $\bar{\mathbf{u}}$, namely,  $\bar{\mathbf{u}} \in \mathbf{L}^{\infty,d/(d-1)}(\Omega)$ and the polynomial degree that is used for its approximation. 
}
\end{remark}

The following result establishes rates of convergence for the errors $\bar{\mathbf{y}} - \bar{\mathbf{y}}_h$, $\bar{p} - \bar{p}_h$, and $\bar{\mathbf{z}} - \bar{\mathbf{z}}_h$.

\begin{theorem}[Rates of convergence]
\label{thm:velocity_rates}
Let $(\bar{\mathbf{y}},\bar{p},\bar{\mathbf{z}},\bar{r},\bar{\mathbf{u}})\in \mathbf{H}_0^1(\Omega)\times L^2(\Omega)/\mathbb{R}\times \mathbf{H}_0^1(\rho,\Omega) \times L^2(\rho,\Omega)/\mathbb{R}\times \mathbb{U}_{ad}$ be the solution to the optimality system \eqref{eq:weak_pde}, \eqref{eq:adj_eq}, and \eqref{eq:variational_ineq} and $(\bar{\mathbf{y}}_h,\bar{p}_h,\bar{\mathbf{z}}_h,\bar{r}_h,\bar{\mathbf{u}}_h)\in \mathbf{V}_h\times Q_h\times\mathbf{V}_h\times Q_h\times \mathbb{U}_{ad}^h$ its numerical approximation given as the solution to \eqref{eq:discrete_state_eq}--\eqref{eq:discrete_adj_eq}. \EO{Then,
\begin{multline}
\label{global_reliability_velocity}
\| \bar{\mathbf{y}} - \bar{\mathbf{y}}_h  \|_{\mathbf{L}^\infty(\Omega)}
\lesssim
h^{2-d/2}\sum_{t\in\mathcal{Z}} |\bar{\mathbf{y}}(t) - \mathbf{y}_{t}| \|\delta_{t}\|_{\mathcal{M}(\Omega)}
\\
+ h^{2-d/2}|\log h| + h|\log h|^{3}\| \bar{\mathbf{u}} \|_{\mathbf{L}^{\vartheta}(\Omega)},
\end{multline}
\vspace{-1em}
\begin{multline}
\label{global_reliability_pressure}
\| \bar{p} - \bar{p}_h  \|_{L^{2}(\Omega)}
\lesssim 
h^{2-d/2}\sum_{t\in\mathcal{Z}} |\bar{\mathbf{y}}(t) - \mathbf{y}_{t}| \|\delta_{t}\|_{\mathcal{M}(\Omega)}
\\
+ h^{2-d/2}|\log h| + h|\log h|^{3}\| \bar{\mathbf{u}} \|_{\mathbf{L}^{\vartheta}(\Omega)},
\end{multline}
and
\begin{multline}\label{global_reliability_adj_vel}
\| \bar{\mathbf{z}} - \bar{\mathbf{z}}_h  \|_{\mathbf{L}^{2}(\Omega)}
\lesssim
h^{2-d/2}\sum_{t\in\mathcal{Z}} |\bar{\mathbf{y}}(t) - \mathbf{y}_{t}| \|\delta_{t}\|_{\mathcal{M}(\Omega)}
\\
+ h^{2-d/2}|\log h| + h|\log h|^{3}\| \bar{\mathbf{u}} \|_{\mathbf{L}^{\vartheta}(\Omega)},
\end{multline}
with $\vartheta > d$. The} hidden constants are independent of the continuous and discrete solutions, the size of the elements in the mesh $\mathscr{T}_h$, and $\#\mathscr{T}_h$. The constants, however, blow up as $\lambda\downarrow 0$.
\end{theorem}

{\it Proof.}
We first control the error $\bar{\mathbf{y}} - \bar{\mathbf{y}}_h$. To accomplish this task, we invoke the pair $(\hat{\mathbf{y}}_{h},\hat p_h)$, defined as the solution to \eqref{eq:state_hat}, and use the triangle inequality to write
\begin{equation*}
\| \bar{\mathbf{y}} - \bar{\mathbf{y}}_h  \|_{\mathbf{L}^\infty(\Omega)} \leq \| \bar{\mathbf{y}} - \hat{\mathbf{y}}_h  \|_{\mathbf{L}^\infty(\Omega)} + \| \hat{\mathbf{y}}_h - \bar{\mathbf{y}}_h \|_{\mathbf{L}^\infty(\Omega)}.
\end{equation*}
The first term on the right hand side of the previous expression is bounded in \eqref{eq:ineq_7}. In view of \eqref{eq:ineq_xitau}, the second term can be bounded by $\|\bar{\mathbf{u}} - \bar{\mathbf{u}}_{h}\|_{\mathbf{L}^{2}(\Omega)}$. The desired estimate \eqref{global_reliability_velocity} follows by collecting the previous estimates and the one obtained in Theorem \ref{thm:control_rates}.

We now control $\bar{p} -\bar{p}_h$. We invoke, again, the auxiliary variable \EO{$(\hat{\mathbf{y}}_h, \hat{p}_h)$, defined as the solution to \eqref{eq:state_hat}, and the triangle inequality to obtain
\begin{equation*}
\| \bar{p} - \bar{p}_h  \|_{L^{2}(\Omega)} \lesssim \| \bar{p} - \hat{p}_h  \|_{L^{2}(\Omega)} + \| \hat{p}_h - \bar{p}_h  \|_{L^{2}(\Omega)}.
\end{equation*}
Since $(\hat{\mathbf{y}}_h,\hat{p}_{h})$ corresponds to the Galerkin approximation of $(\bar{\mathbf{y}},\bar{p})$, then \cite[Proposition 4.16]{Guermond-Ern} yields the error estimate}
\begin{equation*}
\|\bar{p} - \hat{p}_{h} \|_{L^{2}(\Omega)} \lesssim h \left( \|\bar{\mathbf{y}}\|_{\mathbf{H}^{2}(\Omega)} + \|\bar{p}\|_{H^{1}(\Omega)} \right) 
\FF{\lesssim
h\|\bar{\mathbf{u}}\|_{\mathbf{L}^{2}(\Omega)}},
\end{equation*}
\FF{where to obtain the last estimate we have used \eqref{eq:H2regularity}.}
On the other hand, 
\EO{since $\mathbf{V}_h\times Q_h$ satisfy the so--called inf-sup condition \cite[Proposition 4.13]{Guermond-Ern}, we can thus invoke the discrete problem that  $(\hat{\mathbf{y}}_h - \bar{\mathbf{y}}_h,\hat{p}_h - \bar{p}_h)$ solves to obtain
%
\begin{equation*}
\|\hat{p}_h - \bar{p}_h\|_{L^{2}(\Omega)} \lesssim \|\bar{\mathbf{u}} - \bar{\mathbf{u}}_{h}\|_{\mathbf{L}^{2}(\Omega)} + \| \nabla ( \bar{\mathbf{y}}_h - \hat{\mathbf{y}}_h )  \|_{\mathbf{L}^{2}(\Omega)} \lesssim \|\bar{\mathbf{u}} - \bar{\mathbf{u}}_{h}\|_{\mathbf{L}^{2}(\Omega)}.
\end{equation*}
The} result of Theorem \ref{thm:control_rates} allows us to control $\|\bar{\mathbf{u}} - \bar{\mathbf{u}}_{h}\|_{\mathbf{L}^{2}(\Omega)}$.

Finally, we bound $\bar{\mathbf{z}} -\bar{\mathbf{z}}_h$. We begin with the estimate
\begin{equation*}
\| \bar{\mathbf{z}} - \bar{\mathbf{z}}_h  \|_{\mathbf{L}^{2}(\Omega)}
\leq
\| \bar{\mathbf{z}} - \tilde{\mathbf{z}}_h  \|_{\mathbf{L}^{2}(\Omega)} + \| \tilde{\mathbf{z}}_h - \bar{\mathbf{z}}_h  \|_{\mathbf{L}^{2}(\Omega)},
\end{equation*}
where $(\tilde{\mathbf{z}}_h, \tilde r_h)$ solves \eqref{eq:adjoint_tilde}. The term $\| \bar{\mathbf{z}} - \tilde{\mathbf{z}}_h  \|_{\mathbf{L}^{2}(\Omega)}$ can be estimated by using the result of Lemma \ref{lemma:first_estimate}. The remaining term is bounded by using the result of Theorem \ref{thm:weighted_poincare} in conjunction with the stability, in weighted spaces, of the discrete Stokes system \cite[Theorem 4.1]{2019arXiv190500476D}:
\begin{equation*}
\| \tilde{\mathbf{z}}_h - \bar{\mathbf{z}}_h  \|_{\mathbf{L}^{2}(\Omega)}
\lesssim
\| \nabla (\tilde{\mathbf{z}}_h - \bar{\mathbf{z}}_h)  \|_{\mathbf{L}^{2}(\rho,\Omega)}
\lesssim
\| \bar{\mathbf{y}} - \bar{\mathbf{y}}_h\|_{\mathbf{L}^\infty(\Omega)}.
\end{equation*}
The proof concludes by invoking \eqref{global_reliability_velocity}.
\qed


\subsection{\EO{A Semi Discrete Scheme: Error Estimates}}

\EO{In this section, we propose a semidiscrete scheme for the pointwise tracking optimal control problem that is based on the so-called variational discretization approach \cite{MR2122182}. This approach discretizes only the state space (the control space $\mathbb{U}_{ad}$ is not discretized) and induces a discretization of the optimal control variable by projecting the optimal discrete adjoint state into the admissible control set.

The semidiscrete scheme reads as follows: Find min $J(\mathbf{y}_h,\mathbf{q})$ subject to the discrete state equations
\begin{equation}\label{eq:discrete_state_eq_SD}
\left\{
\begin{array}{rcll}
a(\mathbf{y}_h,\mathbf{v}_h) + b(\mathbf{v}_h,p_h) & = & (\mathbf{q},\mathbf{v}_h)_{\mathbf{L}^2(\Omega)} & \quad \forall \mathbf{v}_h \in \mathbf{V}_h, \\
b(\mathbf{y}_h,q_h) & = & 0 & \quad \forall q_h \in Q_h,
\end{array}
\right.
\end{equation}
and the control constraints $\mathbf{q} \in \mathbb{U}_{ad}$. Standard arguments yield the existence and uniqueness of an optimal solution $((\bar{\mathbf{y}}_{h}, \bar{p}_h),\bar{\mathbf{q}})$. In addition,  $((\bar{\mathbf{y}}_{h}, \bar{p}_h),\bar{\mathbf{q}})$ is optimal for the semidiscrete scheme if and only if $(\bar{\mathbf{y}}_{h}, \bar{p}_h)$ solves \eqref{eq:discrete_state_eq_SD}, with $\mathbf{q}$ replaced by $\bar{\mathbf{q}}$, and $\bar{\mathbf{q}}$ satisfies the variational inequality
\begin{equation}\label{eq:discrete_variational_ineq_SD}
(\bar{\mathbf{z}}_h + \lambda \bar{\mathbf{q}},\mathbf{q} - \bar{\mathbf{q}})_{\mathbf{L}^{2}(\Omega)} \geq 0  \quad \forall \mathbf{q} \in \mathbb{U}_{ad},
\end{equation}
where $(\bar{\mathbf{z}}_h, \bar{r}_h)$ solves \eqref{eq:discrete_adj_eq} with ${\mathbf{y}}_h$ replaced by $\bar{\mathbf{y}}_h$.

In the next result we provide error estimates for the semidiscrete scheme.

\begin{theorem}[Rates of convergence for $\bar{\mathbf{u}}$]
\label{thm:control_rates_SD}
Let $(\bar{\mathbf{y}},\bar{p},\bar{\mathbf{z}},\bar{r},\bar{\mathbf{u}})\in \mathbf{H}_0^1(\Omega)\times L^2(\Omega)/\mathbb{R}\\ \times \mathbf{H}_0^1(\rho,\Omega) \times L^2(\rho,\Omega)/\mathbb{R}\times \mathbb{U}_{ad}$ be the solution to the optimality system \eqref{eq:weak_pde}, \eqref{eq:adj_eq}, and \eqref{eq:variational_ineq} and $(\bar{\mathbf{y}}_h,\bar{p}_h,\bar{\mathbf{z}}_h,\bar{r}_h,\bar{\mathbf{q}})\in \mathbf{V}_h\times Q_h\times\mathbf{V}_h\times Q_h\times \mathbb{U}_{ad}$ its numerical approximation given as the solution to the semidiscrete scheme. Then, we have
\begin{multline}\label{eq:global_reliability_2_SD}
\| \bar{\mathbf{u}} - \bar{\mathbf{q}} \|_{\mathbf{L}^2(\Omega)}
\lesssim
h^{2-d/2}\sum_{t\in\mathcal{Z}} |\bar{\mathbf{y}}(t) - \mathbf{y}_{t}|\|\delta_{t}\|_{\mathcal{M}(\Omega)}\\
+ h|\log h|^{3}\left(\|\nabla\bar{\mathbf{y}}\|_{\mathbf{L}^{\infty}(\Omega)} + \|\bar{p}\|_{L^{\infty}(\Omega)}\right).
\end{multline}
The hidden constant is independent of the continuous and discrete solutions, the size of the elements in the mesh $\mathscr{T}_h$, and $\#\mathscr{T}_h$. The constant, however, blows up as $\lambda\downarrow 0$.
\end{theorem}
{\it Proof.} Set $\mathbf{u} = \bar{\mathbf{q}}$ and $\mathbf{q} = \bar{\mathbf{u}}$ in the variational inequalities \eqref{eq:variational_ineq} and \eqref{eq:discrete_variational_ineq_SD}, respectively. Add the obtained inequalities to arrive at
\begin{multline*}
\lambda\|\bar{\mathbf{u}}-\bar{\mathbf{q}}\|_{\mathbf{L}^2(\Omega)}^2
\leq
(\bar{\mathbf{z}}-\bar{\mathbf{z}}_h,\bar{\mathbf{q}}-\bar{\mathbf{u}})_{\mathbf{L}^2(\Omega)} 
\\
= (\bar{\mathbf{z}}-\tilde{\mathbf{z}}_h,\bar{\mathbf{q}}-\bar{\mathbf{u}})_{\mathbf{L}^2(\Omega)}
+
(\tilde{\mathbf{z}}_h-\hat{\mathbf{z}}_h,\bar{\mathbf{q}}-\bar{\mathbf{u}})_{\mathbf{L}^2(\Omega)}
 +
 (\hat{\mathbf{z}}_h-\bar{\mathbf{z}}_h,\bar{\mathbf{q}}-\bar{\mathbf{u}})_{\mathbf{L}^2(\Omega)}.
\end{multline*}
The desired estimate \eqref{eq:global_reliability_2_SD} thus follows from the arguments elaborated in the proof of Theorem \ref{thm:control_rates}. \qed
}
\begin{remark}[variational discretization]
\label{rem:better_control_rates_SD}
\EO{There is an improvement over the fully discrete scheme of section \ref{sec:fully_discrete_scheme}. The error estimate of Theorem \ref{thm:control_rates_SD} reads 
\[
\| \bar{\mathbf{u}} - \bar{\mathbf{q}}  \|_{\mathbf{L}^2(\Omega)} \lesssim h|\log h|^3,
\quad
\| \bar{\mathbf{u}} - \bar{\mathbf{q}}  \|_{\mathbf{L}^2(\Omega)} \lesssim h^{1/2},
\]
for $d=2$ and $d=3$, respectively. Notice that, for $d=2$, this error estimate coincides with the one obtained in Theorem \ref{thm:control_rates}. In particular, it is nearly--optimal in terms of approximation. When $d=3$, \eqref{eq:global_reliability_2_SD} is optimal in terms of regularity and improves upon \eqref{eq:global_reliability_2}.}
\end{remark}


\section{The Optimal Control Problem with Singular Sources}\label{sec:singular_sources_control_probl}

In this section, we precisely describe and analyze the optimal control problem with point sources \eqref{def:cost_func_2}--\eqref{def:box_constraints_2} introduced in Section \ref{sec:singular_source}. We begin by defining the weight $\rho$ as in \eqref{def:weight_rho}--\eqref{def:weight_rho_complete} with obvious modifications that basically entails replacing $\mathcal{E}$ by $\mathcal{D}$. We recall that the cost functional $\mathfrak{J}$ and the set of admissible controls $\mathfrak{U}_{ad}$ are defined by \eqref{def:cost_func_2} and \eqref{def:box_constraints_2}, respectively. \EO{Before we proceed with our analysis, we comment that when deriving a priori error estimates for suitable finite element approximations of problem \eqref{def:cost_func_2}--\eqref{def:box_constraints_2}, it will be essential to assume the existence of $d_\mathcal{D} > 0$ such that
\begin{equation}
\label{eq:assumption_point_sources}
\mathrm{dist} ( \mathcal{D}, \partial \Omega ) \geq d_{\mathcal{D}}.
\end{equation}

The} weak version of the optimal control problem with point sources reads as follows: Find
\begin{equation}\label{def:weak_ocp_ss}
\min \{ \mathfrak{J}(\mathbf{y},\mathcal{U}): (\mathbf{y},\mathcal{U}) \in \mathbf{H}_0^1(\rho,\Omega)\times\mathfrak{U}_{ad} \},
\end{equation}
subject to the following weak formulation of the state equation \eqref{def:state_eq_2}: Find $(\mathbf{y},p) \in \mathbf{H}_0^1(\rho,\Omega) \times  L^2(\rho,\Omega)/\mathbb{R}$ such that
\begin{equation}\label{eq:weak_pde_ss}
\left\{\begin{array}{rcll}
a(\mathbf{y},\mathbf{v})+b(\mathbf{v},p) & = & \displaystyle{\sum_{t\in\mathcal{D}}}\langle \mathbf{u}_{t}\delta_{t},\mathbf{v} \rangle & \quad \forall \mathbf{v} \in \mathbf{H}_0^1(\rho^{-1},\Omega), \\
b(\mathbf{y},q) & = & 0 & \quad \forall q \in L^2(\rho^{-1},\Omega)/\mathbb{R}.
\end{array}
\right.
\end{equation}
We recall that $\langle \cdot,\cdot \rangle$ denotes the duality pairing between $\mathbf{H}_0^1(\rho^{-1},\Omega)'$ and $\mathbf{H}_0^1(\rho^{-1},\Omega)$. If $\alpha \in (d-2,d)$ problem \eqref{eq:weak_pde_ss} is well--posed; see Section \ref{subsub:weak_formulation} for details. Finally, we mention that in view of the continuous embedding of Theorem \ref{thm:weighted_poincare}, $\mathfrak{J}$ is well defined over $\mathbf{H}_0^1(\rho,\Omega)\times\mathfrak{U}_{ad}$. \EO{This further restricts $\alpha$ to belong to the interval} $(d-2,2) \subset (d-2,d)$.

To analyze the optimal control problem with point sources we introduce the so-called control-to-state operator 
\[
 \mathcal{C}: [\mathbb{R}^{d}]^{l} \rightarrow \mathbf{H}_{0}^{1}(\rho,\Omega), \quad [\mathbb{R}^{d}]^{l} \ni \mathcal{U}  \mapsto \mathbf{y}=\mathcal{C}\mathcal{U} \in \mathbf{H}_{0}^{1}(\rho,\Omega),
\]
where $\mathbf{y}=\mathcal{C}\mathcal{U}$ solves problem \eqref{eq:weak_pde_ss}. Since $\alpha \in (d-2,2)$, the map $\mathcal{C}$ is well defined. We can thus define the reduced cost functional
\begin{equation*}
\mathfrak{j}(\mathcal{U}):=\mathfrak{J}(\mathcal{C}\mathcal{U},\mathcal{U})=\frac{1}{2}\|\mathcal{C}\mathcal{U} - \mathbf{y}_{\Omega}\|_{\mathbf{L}^2(\Omega)}^2 + \frac{\lambda}{2}\sum_{t \in \mathcal{D}}| \mathbf{u}_{t} |^{2}.
\end{equation*}
We immediately conclude that $\mathfrak{j}$ is weakly lower semicontinuous and strictly convex $(\lambda > 0)$. This, combined with the fact that $\mathfrak{U}_{ad}$ is compact, allow us to conclude the existence and uniqueness of an optimal control $\bar{\mathcal{U}} \in \mathfrak{U}_{ad}$ and an optimal state $\bar{\mathbf{y}} =\mathcal{C}\bar{\mathcal{U}} \in \mathbf{H}_{0}^{1}(\rho,\Omega)$ that satisfies \eqref{eq:weak_pde_ss} \cite[Theorem 2.14]{Troltzsch}. In addition, the control variable $\bar{\mathcal{U}}$ is optimal for our optimal control problem if and only if \cite[Lemma 2.21]{Troltzsch}
\begin{equation}\label{eq:var_ineq_ss}
\mathfrak{j}'(\bar{\mathcal{U}})({\mathcal{U}} - \bar{\mathcal{U}}) \geq 0\qquad \forall \mathcal{U} \in \mathfrak{U}_{ad}.
\end{equation}
To explore this variational inequality, we introduce the adjoint pair $(\mathbf{z},r) \in \mathbf{H}_{0}^{1}(\Omega) \times L^{2}(\Omega)/\mathbb{R}$, which satisfies
\begin{equation}\label{eq:adj_eq_ss}
\left\{\begin{array}{rcll}
a(\mathbf{z},\mathbf{w})-b(\mathbf{w},r) & = & (\mathbf{y} - \mathbf{y}_{\Omega},\mathbf{w})_{{\mathbf{L}}^2(\Omega)} & \quad \forall \mathbf{w} \in \mathbf{H}_0^1(\Omega), \\
b(\mathbf{z},s) & = & 0 & \quad \forall s \in L^2(\Omega)/\mathbb{R}.
\end{array}
\right.
\end{equation}
Since $\mathbf{y} - \mathbf{y}_{\Omega} \in \mathbf{L}^{2}(\Omega)$, the well--posedness of problem \eqref{eq:adj_eq_ss} is immediate. Moreover, since $\Omega$ is convex, the results of Proposition \ref{pro:H2regularity} yield $(\mathbf{z},r)\in \mathbf{H}^2(\Omega)\times H^1(\Omega)$. This combined with \cite[Proposition 2.3]{2019arXiv190500476D} reveal that the adjoint pair $(\mathbf{z},r) \in \mathbf{H}_{0}^{1}(\rho^{-1},\Omega) \times L^{2}(\rho^{-1},\Omega)$. This result is important because it allows to set $(\mathbf{v},q) = (\mathbf{z},r)$ as a test function in problem \eqref{eq:weak_pde_ss}.

With these ingredients at hand, we proceed to \EO{further explore} optimality conditions for problem \eqref{def:weak_ocp_ss}--\eqref{eq:weak_pde_ss}. Note that, since $(\mathbf{y},p) \in \mathbf{H}_0^1(\rho,\Omega) \times  L^2(\rho,\Omega)/\mathbb{R} \setminus \mathbf{H}_0^1(\Omega) \times  L^2(\Omega)/\mathbb{R}$, we are not \EO{allow} to set $(\mathbf{y},p)$ as a test function in problem \eqref{eq:adj_eq_ss}. We thus have to proceed on the basis of different arguments.
 
\begin{theorem}[Optimality conditions]
\label{thm:opt_cond_ss}
Let $\alpha \in (d-2,2)$. The pair $(\bar{\mathbf{y}},\bar{\mathcal{U}}) \in \mathbf{H}^{1}_{0}(\rho,\Omega) \times \mathfrak{U}_{ad}$ is optimal for   problem \eqref{def:weak_ocp_ss}--\eqref{eq:weak_pde_ss} if and only if $\bar{\mathbf{y}} = \mathcal{C}\bar{\mathcal{U}}$ and the optimal control $\bar{\mathcal{U}}$ satisfies the variational inequality
\begin{equation}\label{eq:var_ineq_ss_2}
\sum_{t \in \mathcal{D}}(\bar{\mathbf{z}}(t) + \lambda\bar{\mathbf{u}}_{t})\cdot(\mathbf{u}_{t} - \bar{\mathbf{u}}_{t} ) \geq 0 \qquad \forall \mathcal{U} = (\mathbf{u}_{1},...,\mathbf{u}_{l}) \in \mathfrak{U}_{ad},
\end{equation}
where $(\bar{\mathbf{z}},\bar{r}) \in \mathbf{H}_{0}^{1}(\Omega) \times L^{2}(\Omega)/\mathbb{R} $ corresponds to the optimal adjoint state, which solves \eqref{eq:adj_eq_ss} with $\mathbf{y}$ replaced by $\bar{\mathbf{y}} = \mathcal{C}\bar{\mathcal{U}}$.
\end{theorem}
{\it Proof.}
A simple computation shows that the variational inequality \eqref{eq:var_ineq_ss} can be rewritten as
\begin{equation}\label{eq:var_ineq2}
\left(\mathcal{C}\bar{\mathcal{U}} - \mathbf{y}_{\Omega},\mathcal{C}({\mathcal{U}} - \bar{\mathcal{U}})\right)_{\mathbf{L}^2(\Omega)} + \lambda\sum_{t \in \mathcal{D}}\bar{\mathbf{u}}_t\cdot(\mathbf{u}_t - \bar{\mathbf{u}}_t) \geq 0
\end{equation}
for all $\mathcal{U} = (\mathbf{u}_1, \ldots, \mathbf{u}_l) \in \mathfrak{U}_{ad}$. In what follows, to simplify the presentation of the material, we let $\mathbf{\mathbf{y}} = \mathcal{C}\mathcal{U}$. Let us concentrate on the first term of the left hand side of the previous expression. To study such a term, we note that $(\mathbf{y} - \bar{\mathbf{y}},p - \bar{p})$ solves
\begin{equation}\label{eq:var_ineq_1}
a(\mathbf{y}-\bar{\mathbf{y}},\mathbf{v})+b(\mathbf{v},p-\bar{p})=\sum_{t \in \mathcal{D}}\langle (\mathbf{u}_{t}-\bar{\mathbf{u}}_t)\delta_{t},\mathbf{v} \rangle, \quad b(\mathbf{y} -\bar{\mathbf{y}},q) = 0
\end{equation}
for all $\mathbf{v}\in \mathbf{H}_0^1(\rho^{-1},\Omega)$ and $q\in L^2(\rho^{-1},\Omega)/\mathbb{R}$, respectively.
Since the variable $\bar{\mathbf{z}} \in \mathbf{H}_{0}^{1}(\rho^{-1},\Omega)$, we are allowed to set $\mathbf{v} = \bar{\mathbf{z}}$ and $q = 0$ in \eqref{eq:var_ineq_1}. This yields
\begin{equation}\label{eq:opt_cond_1}
a(\mathbf{y} - \bar{\mathbf{y}},\bar{\mathbf{z}}) = \sum_{t \in \mathcal{D}}\langle (\mathbf{u}_{t}-\bar{\mathbf{u}}_t)\delta_{t},\bar{\mathbf{z}} \rangle.
\end{equation}
With this identity at hand, a density argument allows us to conclude
\begin{equation}\label{eq:opt_cond_2}
a(\mathbf{y} - \bar{\mathbf{y}},\bar{\mathbf{z}})=(\bar{\mathbf{y}}-\mathbf{y}_\Omega,\mathbf{y}-\bar{\mathbf{y}})_{{\mathbf{L}}^2(\Omega)}.
\end{equation}
In fact, let $\{\mathsf{y}_{n}\}_{n \in \mathbb{N}} \subset \mathbf{C}^{\infty}_{0}(\Omega)$ be such that $\mathsf{y}_n\rightarrow \mathbf{y}-\bar{\mathbf{y}}$ in $\mathbf{H}_0^1(\rho,\Omega)$. Since, for  $n \in \mathbb{N}$, $\mathsf{y}_{n}$ is smooth, we can set $\mathbf{w}= \mathsf{y}_{n}$ and $s=0$ in \eqref{eq:adj_eq_ss}. This yields 
\begin{equation*}
a(\bar{\mathbf{z}},\mathsf{y}_n) - b(\mathsf{y}_n,\bar{r})=(\bar{\mathbf{y}}-\mathbf{y}_\Omega,\mathsf{y}_n)_{{\mathbf{L}}^2(\Omega)}.
\end{equation*}
Now, observe that, on the basis of Theorem \ref{thm:weighted_poincare}, we have
\begin{multline*}
| \left (\bar{\mathbf{y}} - \mathbf{y}_{\Omega},\mathbf{y} - \bar{\mathbf{y}})_{\mathbf{L}^2(\Omega)} - (\bar{\mathbf{y}} - \mathbf{y}_{\Omega},\mathsf{y}_{n})_{{\mathbf{L}}^2(\Omega)} \right|
\\
 \lesssim \|\bar{\mathbf{y}} - \mathbf{y}_{\Omega} \|_{\mathbf{L}^{2}(\Omega)}\|\nabla((\mathbf{y} - \bar{\mathbf{y}}) - \mathsf{y}_{n} )\|_{\mathbf{L}^{2}(\rho,\Omega)} \rightarrow 0, \quad n \uparrow \infty.
\end{multline*}
On the other hand, since $\bar{r} \in L^2(\rho^{-1},\Omega)$, we can  set $q = \bar{r}$ in \eqref{eq:var_ineq_1} to arrive at $b(\mathbf{y} - \bar{\mathbf{y}}, \bar{r}) = 0$. This and the continuity of the bilinear form $b$ on $\mathbf{H}_0^1(\rho,\Omega) \times L^2(\rho^{-1},\Omega)$ imply that $b(\mathsf{y}_n,\bar{r})$ converges to $0$ as $n \uparrow \infty$. Finally, since $\bar{\mathbf{z}}\in \mathbf{H}_0^1(\rho^{-1},\Omega)$, the continuity of bilinear form $a$ on $\mathbf{H}_0^1(\rho,\Omega) \times \mathbf{H}_0^1(\rho^{-1},\Omega)$ allows us to conclude that $a((\mathbf{y}-\bar{\mathbf{y}})-\mathsf{y}_n,\bar{\mathbf{z}})$ tends to $0$ as $n \uparrow \infty$. The collection of these arguments yield the required identity \eqref{eq:opt_cond_2}. 

The proof concludes upon using \eqref{eq:var_ineq2}, \eqref{eq:opt_cond_1}, and \eqref{eq:opt_cond_2}.
\qed

We now introduce, for each $t \in \mathcal{D}$, the projection operator 
\[
 \Pi_{[\mathbf{a}_{t},\mathbf{b}_{t}]}: \mathbb{R}^{d} \rightarrow \mathbb{R}^{d}, \quad \Pi_{[\mathbf{a}_{t},\mathbf{b}_{t}]}(\mathbf{v}) := \min\{\mathbf{b}_{t},\max\{\mathbf{v},\mathbf{a}_{t} \} \}.
\]
With this operator at hand, similar arguments to the ones elaborated in the proof of \cite[Lemma 2.26]{Troltzsch} reveal that $\bar{\mathcal{U}}=(\bar{\mathbf{u}}_1,\ldots,\bar{\mathbf{u}}_{l})$ satisfies \eqref{eq:var_ineq_ss_2} if and only if 
\begin{equation*}
\bar{\mathbf{u}}_{t} = \Pi_{[\mathbf{a}_{t},\mathbf{b}_{t}]}\left(-\lambda^{-1}\bar{\mathbf{z}}(t) \right), \quad t\in \mathcal{D}.
\end{equation*}

To summarize, the pair $(\bar{\mathbf{y}},\bar{\mathcal{U}})$ is optimal for problem \eqref{def:cost_func_2}--\eqref{def:box_constraints_2} if and only if $(\bar{\mathbf{y}},\bar{p},\bar{\mathbf{z}},\bar{r},\bar{\mathcal{U}})\in \mathbf{H}_0^1(\rho,\Omega)\times L^2(\rho,\Omega)/\mathbb{R}\times \mathbf{H}_0^1(\Omega)\times L^2(\Omega)/\mathbb{R}\times \mathfrak{U}_{ad}$ solves \eqref{eq:weak_pde_ss}, \eqref{eq:adj_eq_ss}, and \eqref{eq:var_ineq_ss_2}.

\subsection{Discretization and Error Estimates}
\label{sec:approx_results_2}
We begin by introducing the discrete counterpart of \eqref{def:cost_func_2}--\eqref{def:box_constraints_2}, which reads as follows: Find min $\mathfrak{J}(\mathbf{y}_h,\mathcal{U}_h)$ subject to the discrete state equations
\begin{equation}
\label{eq:discrete_state_eq_ss}
\left\{\begin{array}{rcll}
a(\mathbf{y}_{h},\mathbf{v}_{h})+b(\mathbf{v}_{h},p_{h}) & = & \displaystyle{\sum_{t\in\mathcal{D}}}\langle \mathbf{u}_{h,t}\delta_{t},\mathbf{v}_{h} \rangle & \quad \forall \mathbf{v}_{h} \in \mathbf{V}_h, \\
b(\mathbf{y}_{h},q_{h}) & = & 0 & \quad \forall q_{h} \in Q_{h},
\end{array}
\right.
\end{equation}
and the control constraints $\mathcal{U}_h\in\mathfrak{U}_{ad}$. The spaces $\mathbf{V}_h$ and $Q_{h}$ are given by \eqref{def:discrete_spaces_TH}. \EO{We restrict ourselves to consider such a pair of finite element spaces in order to apply the recent local error estimates of \cite{2019arXiv190706871B}.} We comment that no discretization is needed for the optimal control variable, since the admissible set $\mathfrak{U}_{ad}$ is a subset of a finite dimensional space. 

Standard arguments reveal the existence of a unique optimal pair $(\bar{\mathbf{y}}_h,\bar{\mathcal{U}}_h)$. In addition, the pair $(\bar{\mathbf{y}}_h,\bar{\mathcal{U}}_h)$ is optimal for the previously defined discrete optimal control problem if and only if $\bar{\mathbf{y}}_h$ solves \eqref{eq:discrete_state_eq_ss} and $\bar{\mathcal{U}}_h$ satisfies the variational inequality
\begin{equation}\label{eq:discrete_var_ineq_ss_2}
\sum_{t \in \mathcal{D}}(\bar{\mathbf{z}}_{h}(t) + \lambda\bar{\mathbf{u}}_{h,t})\cdot(\mathbf{u}_{t} - \bar{\mathbf{u}}_{h,t} ) \geq 0 \qquad \forall  \mathcal{U} = (\mathbf{u}_{1},...,\mathbf{u}_{l}) \in \mathfrak{U}_{ad},
\end{equation}
where $(\bar{\mathbf{z}}_h, \bar{r}_h) \in \mathbf{V}_h \times Q_{h}$ solves
\begin{equation}\label{eq:discrete_adj_eq_ss}
\left\{\begin{array}{rcll}
a(\bar{\mathbf{z}}_{h},\mathbf{w}_{h})-b(\mathbf{w}_{h},\bar{r}_{h}) & = & (\bar{\mathbf{y}}_{h} - \mathbf{y}_{\Omega},\mathbf{w}_{h})_{{\mathbf{L}}^2(\Omega)} & \quad \forall  \mathbf{w}_{h} \in \mathbf{V}_{h}, \\
b(\bar{\mathbf{z}}_{h},s_{h}) & = & 0 & \quad \forall  s_{h} \in Q_{h}.
\end{array}
\right.
\end{equation}

To provide an error analysis for the previous scheme we introduce the following \EO{auxiliary} problem: Find $(\hat{\mathbf{y}}_h,\hat{p}_h)\in \mathbf{V}_h\times Q_h$ such that
\begin{align}\label{eq:state_hat_2}
\begin{cases}
\begin{array}{rcll}
a(\hat{\mathbf{y}}_h,\mathbf{v}_h)+b(\mathbf{v}_h,\hat{p}_h)&=&\sum_{t\in\mathcal{D}}\langle \bar{\mathbf{u}}_t\delta_t,\mathbf{v}_h\rangle & \quad \forall  \mathbf{v}_h\in \mathbf{V}_h,\\
b(\hat{\mathbf{y}}_h,q_h) &=&0 & \quad \forall  q_h\in Q_h.
\end{array}
\end{cases}
\end{align}

To simplify the presentation of the material, we define, for $\mathcal{U} =  (\mathbf{u}_{1},...,\mathbf{u}_{l})  \in \mathfrak{U}_{ad}$ and $\mathcal{V} =  (\mathbf{v}_{1},...,\mathbf{v}_{l})  \in \mathfrak{U}_{ad}$, 
\begin{equation*}
\langle \mathcal{U}, \mathcal{V} \rangle_{\mathcal{D}}:= \sum_{t \in \mathcal{D}} \mathbf{u}_t \cdot \mathbf{v}_t,
\qquad 
\VERT \mathcal{U} \VERT_{\mathcal{D}}:= \sqrt{\langle \mathcal{U}, \mathcal{U} \rangle } = \left( \sum_{t\in\mathcal{D}} | \mathbf{u}_t |^2 \right)^{\frac{1}{2}}.
\end{equation*}
If $\mathbf{w} \in \mathbf{C}(\bar{\Omega})$ and $\mathcal{V} = (\mathbf{v}_{1},...,\mathbf{v}_{l})  \in \FF{[\mathbb{R}^{d}]^{l}}$, $\langle \mathbf{w}, \mathcal{V} \rangle_{\mathcal{D}} :=  \sum_{t \in \mathcal{D}} \mathbf{w}(t) \cdot \mathbf{v}_t.$

With the discrete system \eqref{eq:discrete_state_eq_ss}--\eqref{eq:discrete_adj_eq_ss} at hand, we are in conditions to present the main result of this section. \EO{Our arguments are inspired by ideas developed in \cite{MR3225501,MR3116646} and yield a nearly--optimal error estimate in terms of approximation.}

\begin{theorem}[Rates of convergence for $\mathcal{\bar{U}}$]
\label{thm:control_rates_2}
Let $(\bar{\mathbf{y}},\bar{p},\bar{\mathbf{z}},\bar{r},\bar{\mathcal{U}})\in \mathbf{H}_0^1(\rho,\Omega) \times L^2(\rho,\Omega)/\mathbb{R}\times \mathbf{H}_0^1(\Omega) \times L^2(\Omega)/\mathbb{R}\times \mathfrak{U}_{ad}$ be the solution to the optimality system \eqref{eq:weak_pde_ss}, \eqref{eq:adj_eq_ss}, and \eqref{eq:var_ineq_ss_2} and $(\bar{\mathbf{y}}_h,\bar{p}_h,\bar{\mathbf{z}}_h,\bar{r}_h,\bar{\mathcal{U}}_h)\in \mathbf{V}_h\times Q_h\times\mathbf{V}_h\times Q_h\times \mathfrak{U}_{ad}$ its numerical approximation given as the solution to \eqref{eq:discrete_state_eq_ss}--\eqref{eq:discrete_adj_eq_ss}. \EO{Let $\Omega_0$, $\Omega_1 \subset \Omega$ be such that $\mathcal{D} \Subset \Omega_0 \subset \Omega_1$ with $\mathrm{dist}(\bar \Omega_0,\partial \Omega_1) \geq \mathfrak{d} \geq \mathfrak{c}h$ for $\mathfrak{c}$ sufficiently large. If $d=2$ and $\mathbf{y}_\Omega\in \mathbf{L}^{\kappa}(\Omega)$, for every $\kappa \in (2,\infty)$, then 
\begin{equation*}
\VERT \bar{\mathcal{U}} - \bar{\mathcal{U}}_h \VERT_{\mathcal{D}}
\lesssim
 |\log h|^3 h^{2}  \left( \| \mathbf{y}_{\Omega}\|_{\mathbf{L}^{\kappa}(\Omega)} + \sum_{t\in\mathcal{D}} 
 |\bar{\mathbf{u}}_t| \|\delta_{t}\|_{\mathcal{M}(\Omega)} + \VERT \bar{\mathcal{U}} \VERT_{\mathcal{D}} \right).
\end{equation*}
The} hidden constant is independent of the continuous and discrete solutions, the size of the elements in the mesh $\mathscr{T}_h$, and $\#\mathscr{T}_h$. The constant, however, blows up as $\lambda\downarrow 0$.
\end{theorem}
{\it Proof.} We proceed in $4$ steps.

\underline{Step 1.} Set $\mathcal{U}=\bar{\mathcal{U}}_h$ in \eqref{eq:var_ineq_ss_2} and $\mathcal{U}=\bar{\mathcal{U}}$ in \eqref{eq:discrete_var_ineq_ss_2}. Adding the obtained inequalities we arrive at the basic estimate
\begin{equation}
\begin{aligned}
\lambda \VERT \bar{\mathcal{U}} - \bar{\mathcal{U}}_h \VERT^2_{\mathcal{D}} & = 
\lambda\sum_{t\in\mathcal{D}}|\bar{\mathbf{u}}_{h,t}-\bar{\mathbf{u}}_t|^2
\\
& \leq 
\sum_{t\in\mathcal{D}} \left((\bar{\mathbf{z}}-\bar{\mathbf{z}}_h)(t)\right)\cdot \left( \bar{\mathbf{u}}_{h,t}-\bar{\mathbf{u}}_t\right)
= \langle \bar{\mathbf{z}}-\bar{\mathbf{z}}_h, \bar{\mathcal{U}}_{h}-\bar{\mathcal{U}}\rangle_{\mathcal{D}}.
\end{aligned}
\label{eq:ineq_2_prob_2}
\end{equation}

\underline{Step 2.}
Define the discrete auxiliary variables $(\hat{\mathbf{z}}_h, \hat r_h)$ and $(\tilde{\mathbf{z}}_h, \tilde{r}_h)$ as the solutions to 
\begin{equation*}
\left\{\begin{array}{rcll}
a(\mathbf{w}_h,\hat{\mathbf{z}}_h) - b(\mathbf{w}_h,\hat{r}_h)&=&(\hat{\mathbf{y}}_h-\mathbf{y}_\Omega,\mathbf{w}_h)_{\mathbf{L}^2(\Omega)} &  \quad \forall \mathbf{w}_h \in \mathbf{V}_h,\\
b(\hat{\mathbf{z}}_h,s_h) &=&0 \quad & \quad \forall s_h \in Q_h,
\end{array}
\right.\hspace{-0.45cm}
\end{equation*}
and
\begin{equation*}
\left\{\begin{array}{rcll}
a(\mathbf{w}_h,\tilde{\mathbf{z}}_h) - b(\mathbf{w}_h,\tilde{r}_h)&=&(\bar{\mathbf{y}}-\mathbf{y}_\Omega,\mathbf{w}_h)_{\mathbf{L}^2(\Omega)}& \quad \forall  \mathbf{w}_h \in \mathbf{V}_h,\\
b(\tilde{\mathbf{z}}_h,s_h) &=&0  & \quad \forall  s_h \in Q_h,
\end{array}
\right.\hspace{-0.45cm}
\end{equation*}
respectively, where $(\hat{\mathbf{y}}_h, \hat p_h)$ solves \eqref{eq:state_hat_2}.

We now invoke \eqref{eq:ineq_2_prob_2} and add and subtract $\hat{\mathbf{z}}_h$ and $\tilde{\mathbf{z}}_h$ to obtain 
\begin{multline}\label{eq:ineq_3_prob_2}
\lambda \VERT \bar{\mathcal{U}} - \bar{\mathcal{U}}_h \VERT^2_{\mathcal{D}}
\leq 
\langle \bar{\mathbf{z}}-\tilde{\mathbf{z}}_h, \bar{\mathcal{U}}_{h}-\bar{\mathcal{U}}\rangle_{\mathcal{D}} 
\\
+
\langle \tilde{\mathbf{z}}_h-\hat{\mathbf{z}}_h, \bar{\mathcal{U}}_{h}-\bar{\mathcal{U}}\rangle_{\mathcal{D}}+\langle \hat{\mathbf{z}}_h-\bar{\mathbf{z}}_h, \bar{\mathcal{U}}_{h}-\bar{\mathcal{U}}\rangle_{\mathcal{D}}=:\mathbf{I}+\mathbf{II}+\mathbf{III}.
\end{multline}

Similar arguments to those elaborated in Step 2 of Theorem \ref{thm:control_rates} allow us to conclude that $\mathbf{III} = -\|\bar{\mathbf{y}}_h-\hat{\mathbf{y}}_h\|_{\mathbf{L}^2(\Omega)}^2\leq 0$. As a result, we obtain 
\begin{equation}\label{eq:estimate_I_and_II}
\lambda \VERT \bar{\mathcal{U}} - \bar{\mathcal{U}}_h \VERT^2_{\mathcal{D}}
\leq
\mathbf{I}+\mathbf{II}.
\end{equation}

\underline{Step 3.} \EO{We now estimate the term $\mathbf{I}$. To accomplish this task, we first note that, since $d=2$, the velocity field that solves \eqref{eq:weak_pde_ss} satisfies $\bar{\mathbf{y}}\in \mathbf{W}^{1,\nu}(\Omega)$ for every $\nu<2$. A standard Sobolev embedding result thus implies that $\bar{\mathbf{y}}\in \mathbf{L}^{\kappa}(\Omega)$ for every $\kappa<\infty$. Consequently, $\bar{\mathbf{y}} - \mathbf{y}_{\Omega} \in \mathbf{L}^{\kappa}(\Omega)$ for every $\kappa < \infty$. We can thus apply the interior regularity results of \cite[Theorem IV 4.1]{Gal11} to conclude that $(\bar{\mathbf{z}},\bar{r})\in \mathbf{W}^{2,\kappa}(\Omega_1) \times W^{1,\kappa}(\Omega_1)$ for every  $\kappa < \infty$ together with the bound
\[
| \bar{\mathbf{z}} |_{\mathbf{W}^{2,\kappa}(\Omega_1) }
+
| \bar{r} |_{W^{1,\kappa}(\Omega_1) } \leq C_{\kappa} \| \bar{\mathbf{y}} - \mathbf{y}_{\Omega}\|_{\mathbf{L}^{\kappa}(\Omega)},
\]
where we have also used the first item in Section 5.5 of \cite{MR2321139}. One of the main ingredients in the proof of \cite[Theorem IV 4.1]{Gal11} are the global regularity estimates of \cite[Theorem IV 2.1]{Gal11}, which in turn follow from the Calder\'on--Zygmund theorem. The constant $C_{\kappa}$ can be traced and behaves as the constant involved in Calder\'on--Zygmund theorem: $C_{\kappa} \leq C \kappa$, with $C$ being independent of $\kappa$; see \cite[Remark II.11.2]{Gal11}. On the other hand, since $\Omega$ is convex we have at hand the global regularity result $(\bar{\mathbf{z}}, \bar{r}) \in \mathbf{H}^2(\Omega) \times H^1(\Omega)$; see Proposition \ref{pro:H2regularity}. We are thus in position to apply \cite[Theorem 6.3]{2019arXiv190706871B} combined with \cite[Remark 2.18]{2019arXiv190706871B} to arrive at
\begin{multline}\label{eq:ineq_4_prob_2}
\| \bar{\mathbf{z}} - \tilde{\mathbf{z}}_h \|_{\mathbf{L}^{\infty}(\Omega_0)}
\lesssim \inf \left[
|\log h|\left( |\log h|\|\bar{\mathbf{z}} - \mathbf{w}_h \|_{\mathbf{L}^{\infty}(\Omega_1)} + h\|\bar{r} - r_h\|_{L^{\infty}(\Omega_1)} \right) \right.
\\
\left.
+ |\log h| \left( h \|\bar{\mathbf{z}} - \mathbf{w}_h \|_{\mathbf{H}^1(\Omega)} + \|\bar{\mathbf{z}} - \mathbf{w}_h \|_{\mathbf{L}^2(\Omega)}+ h\|\bar{r} - r_h\|_{L^{2}(\Omega)} \right) \right],
\end{multline}
where the infimum is taken over the whole space $\mathbf{V}_h \times Q_h$. Utilize the aforementioned regularity results and standard interpolation error estimates to obtain
\begin{multline*}
\| \bar{\mathbf{z}} - \tilde{\mathbf{z}}_h \|_{\mathbf{L}^{\infty}(\Omega_0)}
\leq 
C \kappa |\log h|^2 h^{2-2/\kappa}  \| \bar{\mathbf{y}} - \mathbf{y}_{\Omega}\|_{\mathbf{L}^{\kappa}(\Omega)},
\\
+ C |\log h| h^2 \left( \|\bar{\mathbf{z}}\|_{\mathbf{H}^2(\Omega)} + \|\bar{r} \|_{H^{1}(\Omega)} \right),
\end{multline*}
for any $\kappa < \infty$. Inspired by \cite[page 3]{MR637283}, we thus set $\kappa = |\log h|$ to conclude
\begin{equation}\label{eq:ineq_4_prob_2_xx}
\| \bar{\mathbf{z}} - \tilde{\mathbf{z}}_h \|_{\mathbf{L}^{\infty}(\Omega_0)}
\lesssim
 |\log h|^3 h^{2}  \left( \| \mathbf{y}_{\Omega}\|_{\mathbf{L}^{\kappa}(\Omega)} + \sum_{t\in\mathcal{D}} 
 |\bar{\mathbf{u}}_t| \|\delta_{t}\|_{\mathcal{M}(\Omega)} \right),
\end{equation}
where we have also used a stability estimate for the problem that $(\bar{\mathbf{y}},\bar{p})$ solves and the regularity results of Proposition \ref{pro:H2regularity}.  

\underline{Step 4.} We conclude} by estimating the term $\mathbf{II}$ in \eqref{eq:ineq_3_prob_2}. To accomplish this task, we proceed on the basis of a duality argument. Let us define the pair $(\boldsymbol\varphi,\pi)\in \mathbf{H}_0^1(\Omega)\times L^2(\Omega)/\mathbb{R}$ as the solution to
\begin{equation}\label{eq:aux_eq_sgn}
\left\{\begin{array}{rcll}
a(\boldsymbol\varphi,\mathbf{v})+b(\mathbf{v},\pi) & = & (\text{sgn}(\bar{\mathbf{y}}-\hat{\mathbf{y}}_h),\mathbf{v})_{\mathbf{L}^2(\Omega)} &\quad  \forall  \mathbf{v} \in \mathbf{H}_0^1(\Omega), \\
b(\boldsymbol\varphi,q) & = & 0 & \quad \forall q \in L^2(\Omega)/\mathbb{R},
\end{array}
\right.
\end{equation}
where the pair $(\hat{\mathbf{y}}_h,\hat{p}_h)$ solves \eqref{eq:state_hat_2}.  Since $\|\text{sgn}(\bar{\mathbf{y}}-\hat{\mathbf{y}}_h)\|_{\mathbf{L}^\infty(\Omega)}\leq 1$, we can apply, \EO{again, the results of \cite[Theorem IV 4.1]{Gal11} to conclude that $(\boldsymbol\varphi,\pi) \in \mathbf{W}^{2,\infty}(\Omega_1) \times W^{1,\infty}(\Omega_1)$.
We are thus in position to invoke, again, \cite[Theorem 6.3]{2019arXiv190706871B} combined with \cite[Remark 2.18]{2019arXiv190706871B} to conclude that
\begin{multline} 
\|\boldsymbol\varphi-\boldsymbol\varphi_h\|_{\mathbf{L}^\infty(\Omega_0)}
\lesssim |\log h|^2 h^2 \left( \|\boldsymbol\varphi\|_{\mathbf{W}^{2,\infty}(\Omega_1)} + \|r\|_{W^{1,\infty}(\Omega_1)}\right)
\\
+  |\log h| h^2 \left( \|\boldsymbol\varphi\|_{\mathbf{H}^{2}(\Omega)} + \|r\|_{H^{1}(\Omega)}\right)
\lesssim  |\log h|^2 h^2.
\label{eq:eq:ineq_5_prob_2} 
\end{multline}
Here,  $(\boldsymbol \varphi_h, \pi_h)$ corresponds to the finite element approximation of $(\boldsymbol \varphi, \pi)$. Notice that, since $\Omega$ is convex, we immediately have that $(\boldsymbol\varphi, \pi) \in \mathbf{H}^2(\Omega) \times H^1(\Omega)$. To obtain the last inequality, we have used that $\|\text{sgn}(\bar{\mathbf{y}}-\hat{\mathbf{y}}_h)\|_{\mathbf{L}^\infty(\Omega)}\leq 1$.

On} the other hand, we can apply \cite[Proposition 2.3]{2019arXiv190500476D} to obtain that $(\boldsymbol\varphi, \pi) \in \mathbf{H}^2(\Omega)\cap \mathbf{H}_0^1(\Omega) \times H^1(\Omega) \cap L^2(\Omega)/ \mathbb{R} \hookrightarrow \mathbf{H}_0^1(\rho^{-1},\Omega) \times L^2(\rho^{-1},\Omega)$. 
\EO{Since} the pair $(\boldsymbol\varphi, \pi)$ solves \eqref{eq:aux_eq_sgn}, similar arguments to the ones that allowed us to derive \eqref{eq:opt_cond_2} yield
\begin{equation*}
 \|\bar{\mathbf{y}}-\hat{\mathbf{y}}_h\|_{\mathbf{L}^1(\Omega)} =\int_\Omega \text{sgn}(\bar{\mathbf{y}}-\hat{\mathbf{y}}_h)(\bar{\mathbf{y}}-\hat{\mathbf{y}}_h) 
=a(\boldsymbol \varphi,\bar{\mathbf{y}}-\hat{\mathbf{y}}_h)+b(\bar{\mathbf{y}}-\hat{\mathbf{y}}_h,\pi).
\end{equation*}
Note that $(\hat{\mathbf{y}}_h, \hat p_h)$  corresponds to the finite element approximation, within the space \EO{$\mathbf{V}_h \times Q_h$, of $(\bar{\mathbf{y}}, \bar{p})$, the solution to problem \eqref{eq:weak_pde_ss}.} We thus invoke Galerkin orthogonality, twice, and set $\mathbf{v} = \boldsymbol\varphi-\boldsymbol\varphi_h \in  \mathbf{H}_0^1(\rho^{-1},\Omega)$ in \eqref{eq:weak_pde_ss} to arrive at
\begin{equation*}
\|\bar{\mathbf{y}}-\hat{\mathbf{y}}_h\|_{\mathbf{L}^1(\Omega)}
=a(\bar{\mathbf{y}},\boldsymbol\varphi-\boldsymbol\varphi_h)+b(\boldsymbol\varphi-\boldsymbol\varphi_h,\bar{p})
=\langle \boldsymbol\varphi-\boldsymbol\varphi_h, \bar{\mathcal{U}}\rangle_{\mathcal{D}}.
\end{equation*}
Finally, apply \eqref{eq:eq:ineq_5_prob_2} to obtain the error \EO{estimate
\begin{equation}
\|\bar{\mathbf{y}}-\hat{\mathbf{y}}_h\|_{\mathbf{L}^1(\Omega)} 
\lesssim  \|\boldsymbol\varphi-\boldsymbol\varphi_h\|_{\mathbf{L}^\infty(\Omega_0)} 
\VERT \bar{\mathcal{U}} \VERT_{\mathcal{D}}
\lesssim
 |\log h|^2 h^2
\VERT \bar{\mathcal{U}} \VERT_{\mathcal{D}}.
\label{eq:aux_aux_aux}
\end{equation}

With} the previous estimates at hand, we can thus bound $\|\tilde{\mathbf{z}}_h-\hat{\mathbf{z}}_h\|_{\mathbf{L}^\infty(\Omega)}$. To accomplish this task, we invoke a standard inverse estimate \cite[Lemma 4.9.2]{brenner}, the problem that $\tilde{\mathbf{z}}_h-\hat{\mathbf{z}}_h$ solves, and estimate \eqref{eq:aux_aux_aux}. In fact, we have
\begin{align}\label{eq:ineq_6_prob_2}
\begin{split}
\|\tilde{\mathbf{z}}_h-\hat{\mathbf{z}}_h\|_{\mathbf{L}^\infty(\Omega)}^2
&\lesssim
(1+|\log h|)\|\nabla (\tilde{\mathbf{z}}_h-\hat{\mathbf{z}}_h)\|_{\mathbf{L}^2(\Omega)}^2 \\
&\lesssim
(1+|\log h|)\|\bar{\mathbf{y}}-\hat{\mathbf{y}}_h\|_{\mathbf{L}^1(\Omega)}\|\tilde{\mathbf{z}}_h-\hat{\mathbf{z}}_h\|_{\mathbf{L}^\infty(\Omega)}\\
& \lesssim
\EO{ |\log h|^2 h^2(1+|\log h|)
\VERT \bar{\mathcal{U}} \VERT_{\mathcal{D}}}
\|\tilde{\mathbf{z}}_h-\hat{\mathbf{z}}_h\|_{\mathbf{L}^\infty(\Omega)}.
\end{split}
\end{align}

The proof concludes by gathering the estimates \eqref{eq:estimate_I_and_II}, \eqref{eq:ineq_4_prob_2_xx}, and \eqref{eq:ineq_6_prob_2}.
\qed

\begin{theorem}[Rates of convergence]
\label{thm:velocity_rates_2}
Let $(\bar{\mathbf{y}},\bar{p},\bar{\mathbf{z}},\bar{r},\bar{\mathcal{U}})$ be the solution to the optimality system \eqref{eq:weak_pde_ss}, \eqref{eq:adj_eq_ss}, and \eqref{eq:var_ineq_ss_2}, in $\mathbf{H}_0^1(\rho,\Omega) \times L^2(\rho,\Omega)/\mathbb{R}\times \mathbf{H}_0^1(\Omega) \times L^2(\Omega)/\mathbb{R}\times \mathfrak{U}_{ad}$, and $(\bar{\mathbf{y}}_h,\bar{p}_h,\bar{\mathbf{z}}_h,\bar{r}_h,\bar{\mathcal{U}}_h) \in \mathbf{V}_h \times Q_h\times\mathbf{V}_h\times Q_h\times \mathfrak{U}_{ad}$ its numerical approximation given as the solution to \eqref{eq:discrete_state_eq_ss}--\eqref{eq:discrete_adj_eq_ss}. \EO{Let $\Omega_0$, $\Omega_1 \subset \Omega$ be such that $\mathcal{D} \Subset \Omega_0 \subset \Omega_1$ with $\mathrm{dist}(\bar \Omega_0,\partial \Omega_1) \geq \mathfrak{d} \geq \mathfrak{c}h$ for $\mathfrak{c}$ sufficiently large. If $d=2$ and $\mathbf{y}_\Omega \in \mathbf{L}^{\kappa}(\Omega)$ for every $\kappa \in (2,\infty)$, then 
\begin{multline*}
\| \bar{\mathbf{y}} - \bar{\mathbf{y}}_{h}  \|_{\mathbf{L}^{2}(\Omega)}
\lesssim  
(h+ |\log h|^3 h^{2})\sum_{t\in\mathcal{D}} |\bar{\mathbf{u}}_{t}|\|\delta_{t}\|_{\mathcal{M}(\Omega)}
\\
+ |\log h|^3 h^{2}  \left( \| \mathbf{y}_{\Omega}\|_{\mathbf{L}^{\kappa}(\Omega)} + \VERT \bar{\mathcal{U}} \VERT_{\mathcal{D}} \right),
\end{multline*}
and
\begin{multline*}
\| \nabla (\bar{\mathbf{z}} - \bar{\mathbf{z}}_{h}) \|_{\mathbf{L}^{2}(\Omega)}
\lesssim  
h \left( \sum_{t\in\mathcal{D}} |\bar{\mathbf{u}}_{t}|\|\delta_{t}\|_{\mathcal{M}(\Omega)} 
+ 
\| \mathbf{y}_{\Omega}\|_{\mathbf{L}^{2}(\Omega)} \right)
\\
+ |\log h|^3 h^{2}  \left( \| \mathbf{y}_{\Omega}\|_{\mathbf{L}^{\kappa}(\Omega)} 
 + \VERT \bar{\mathcal{U}} \VERT_{\mathcal{D}} \right),
\end{multline*}
with a similar estimate for $\| \bar{r} - \bar{r}_{h} \|_{L^{2}(\Omega)}$. The} hidden constants are independent of the continuous and discrete solutions, the size of the elements in the mesh $\mathscr{T}_h$, and $\#\mathscr{T}_h$. The constants, however, blow up as $\lambda\downarrow 0$.
\end{theorem}
{\it Proof.} We first estimate $\bar{\mathbf{y}} - \bar{\mathbf{y}}_{h}$. To accomplish this task, \EO{we invoke the pair $(\hat{\mathbf{y}}_h, \hat{p}_h)$, defined as the solution to \eqref{eq:state_hat_2}, and write
\begin{equation}\label{eq:des1_vel_rates_2}
\| \bar{\mathbf{y}} - \bar{\mathbf{y}}_{h} \|_{\mathbf{L}^{2}(\Omega)} \lesssim \| \bar{\mathbf{y}} - \hat{\mathbf{y}}_h \|_{\mathbf{L}^{2}(\Omega)} + \| \hat{\mathbf{y}}_h - \bar{\mathbf{y}}_{h} \|_{\mathbf{L}^{2}(\Omega)}.
\end{equation}
Apply now the Poincar\'e inequality of Theorem \ref{thm:weighted_poincare} and the stability of the discrete Stokes system in weighted spaces \cite[Theorem 4.1]{2019arXiv190500476D} to obtain
\[ 
\| \hat{\mathbf{y}}_h - \bar{\mathbf{y}}_{h} \|_{\mathbf{L}^{2}(\Omega)} \lesssim \| \nabla (\hat{\mathbf{y}}_h - \bar{\mathbf{y}}_{h} ) \|_{\mathbf{L}^{2}(\rho,\Omega)} \lesssim 
\VERT \bar{\mathcal{U}} - \bar{\mathcal{U}}_h \VERT_{\mathcal{D}}.
\]
On} the other hand, since $(\hat{\mathbf{y}}_{h},\hat{p}_h)$ corresponds to the Galerkin approximation of $(\bar{\mathbf{y}},\bar{p})$, we estimate the second term on the right hand side of \eqref{eq:des1_vel_rates_2} in view of the error estimate \eqref{eq:apriori_error_delta}. \EO{Collect the derived estimate to arrive at the desired estimate for $\| \bar{\mathbf{y}} - \bar{\mathbf{y}}_{h} \|_{\mathbf{L}^{2}(\Omega)}$.}

Similar arguments can be used to estimate the terms $\| \nabla (\bar{\mathbf{z}} - \bar{\mathbf{z}}_{h}) \|_{\mathbf{L}^{2}(\Omega)}$ and $\| \bar{r} - \bar{r}_{h} \|_{L^{2}(\Omega)}$. This concludes the proof.
\qed

\begin{remark}[Rates of convergence for $\bar{\mathcal{U}}$ and $ \bar{\mathbf{y}}$]
\label{rem:rates_dirac_sources}
\EO{The error estimate of Theorem \ref{thm:control_rates_2} reads 
\[
\VERT \bar{\mathcal{U}} - \bar{\mathcal{U}}_h \VERT_{\mathcal{D}}
\lesssim
 |\log h|^3 h^{2}.
\]
This estimate is nearly--optimal in terms of approximation (nearly because of the presence of the $\log$-term) and improves the one derived in \cite[Theorem 5.1]{MR3800041}, for the Poisson problem, which behaves as $\mathcal{O}(h^{2-\epsilon})$, for every $\epsilon >0$. The error estimate obtained in Theorem \ref{thm:velocity_rates_2} for the discretization of the state velocity field reads as follows: $\| \bar{\mathbf{y}} - \bar{\mathbf{y}}_{h}  \|_{\mathbf{L}^{2}(\Omega)} \lesssim h$. We notice that this estimate is optimal in terms of regularity.
}
\end{remark}

\section{Numerical Examples}\label{sec:numerical_ex}

In this section, we conduct a series of numerical examples that illustrate the performance of the discrete schemes \eqref{eq:discrete_state_eq}--\eqref{eq:discrete_adj_eq} and \eqref{eq:discrete_state_eq_ss}--\eqref{eq:discrete_adj_eq_ss} when approximating the solutions to the optimization problems described in Sections \ref{sec:pointwise_tracking} and \ref{sec:singular_sources_control_probl}, respectively.

\subsection{Implementation}\label{sec:implementation}

All the experiments have been carried out with the help of a code that we implemented using \texttt{C++}. All matrices have been assembled exactly. The right hand sides as well as the approximation errors are computed with the help of a quadrature formula that is exact for polynomials of degree $19$ for two dimensional domains and degree $14$ for three dimensional domains.  
The global linear systems were solved using the multifrontal massively parallel sparse direct solver (MUMPS) \cite{MUMPS1,MUMPS2}. 

\normalsize
\begin{figure}[ht]
\centering
\begin{minipage}{0.315\textwidth}\centering
\includegraphics[trim={0 0 0 0},clip,width=3.0cm,height=2.5cm,scale=0.2]{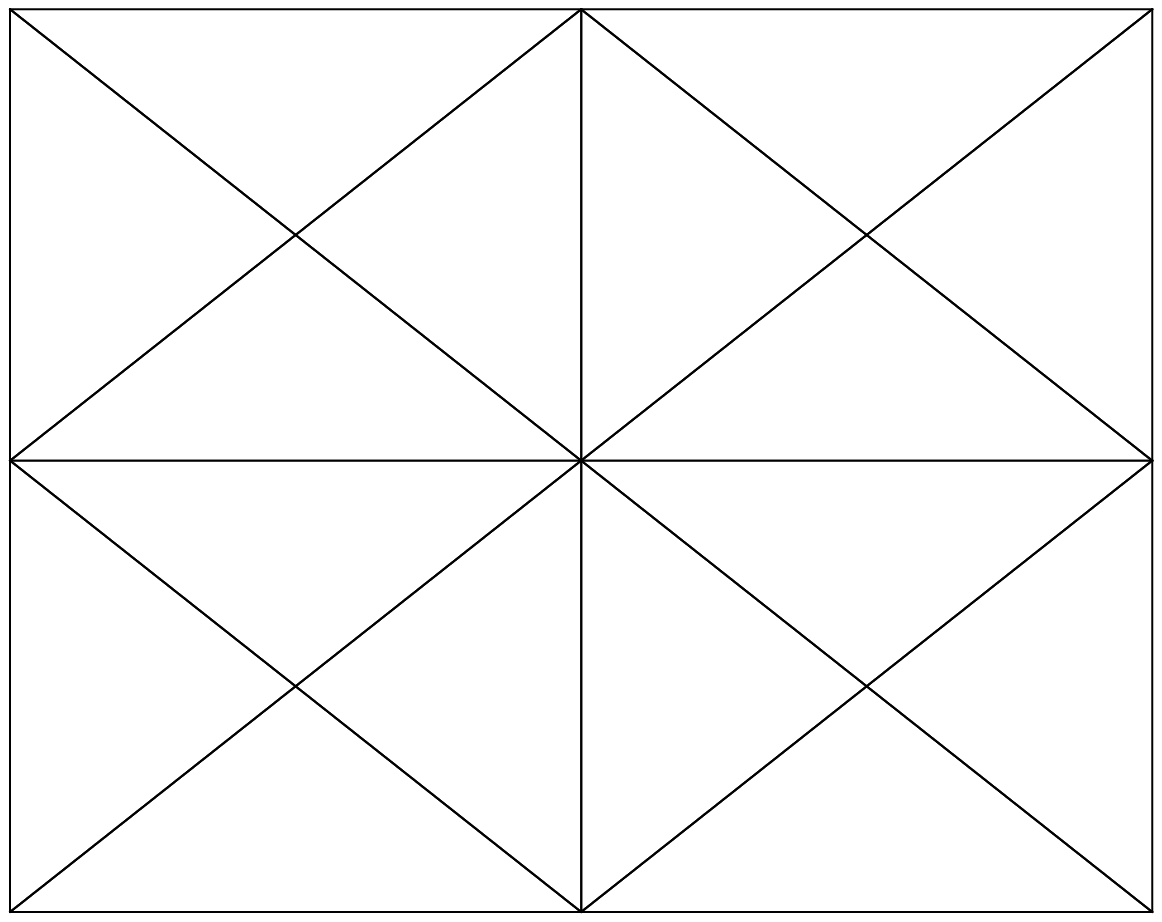}\\
\end{minipage}
\begin{minipage}{0.315\textwidth}\centering
\includegraphics[trim={0 0 0 0},clip,width=2.7cm,height=2.7cm,scale=0.2]{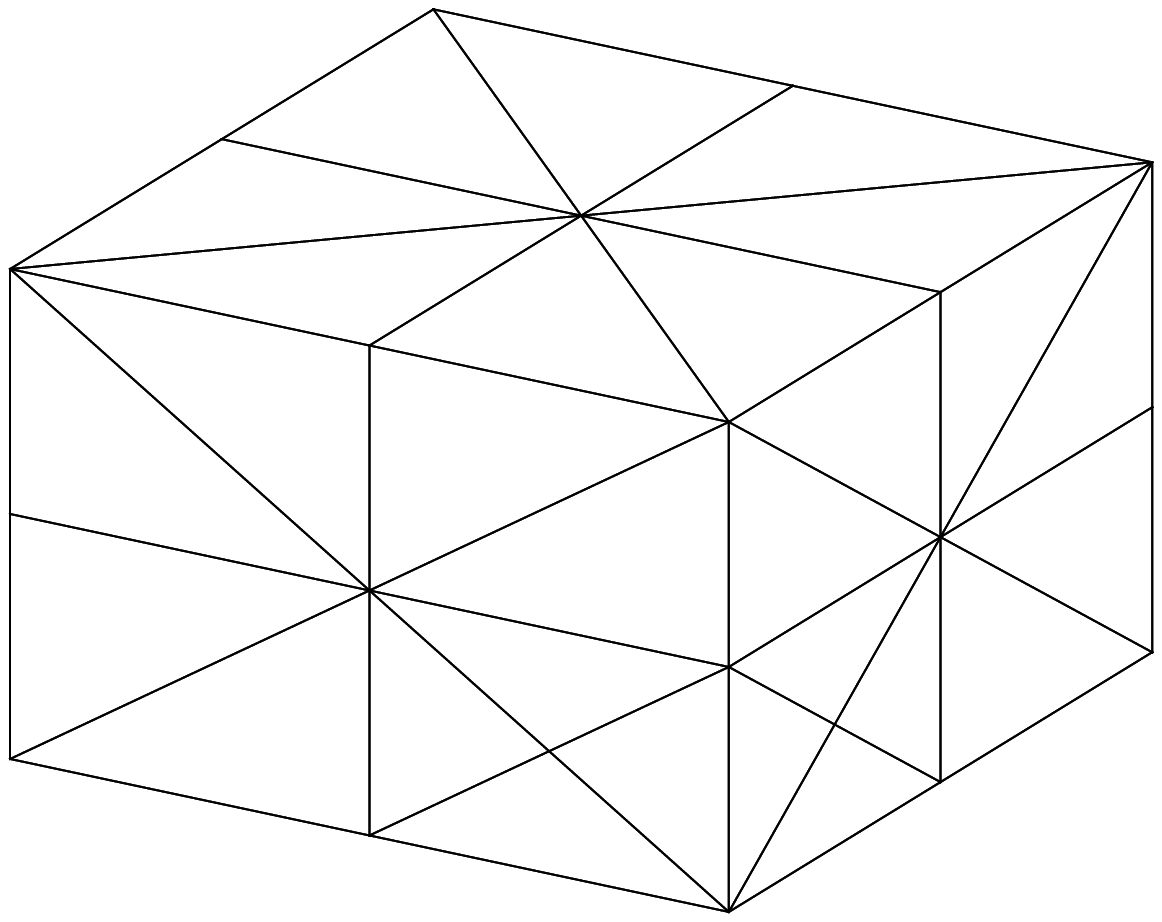}\\
\end{minipage}
\caption{The initial meshes used when the domain $\Omega$ is a square (left) or a cube (right).}
\label{fig:initial_meshes}
\end{figure}

In all the examples we set $\lambda = 1$ and $\Omega = (0,1)^{d}$ with $d\in\{2,3\}$. For a given partition $\mathscr{T}_h$ of $\Omega$, for the first problem, we seek $(\bar{\mathbf{y}}_h,\bar{p}_h,\bar{\mathbf{z}}_h,\bar{r}_h,\bar{\mathbf{u}}^{}_h)\,\in \mathbf{V}_h\times Q_h\times\mathbf{V}_h\times Q_h\times\mathbb{U}_{ad}^h$ that solves the discrete optimality system \eqref{eq:discrete_state_eq}--\eqref{eq:discrete_adj_eq}, while for the second problem, we seek  $(\bar{\mathbf{y}}_h,\bar{p}_h,\bar{\mathbf{z}}_h,\bar{r}_h,\bar{\mathcal{U}}^{}_h)\,\in \mathbf{V}_h\times Q_h\times\mathbf{V}_h\times Q_h\times\mathfrak{U}_{ad}$ that solves the discrete optimality system \eqref{eq:discrete_state_eq_ss}--\eqref{eq:discrete_adj_eq_ss}. In all the numerical examples we make use of the Taylor--Hood element defined in \eqref{def:discrete_spaces_TH}. To solve the associated minimization problems, we use the Newton--type primal--dual active set strategy as described in \cite[Section 2.12.4]{Troltzsch}.

We consider problems with inhomogeneous Dirichlet boundary conditions whose exact solutions are known. Note that this violates the assumption of homogeneous Dirichlet boundary conditions which is needed for the analysis that we have performed, but it retains its essential difficulties and
singularities and allows us to evaluate the rates of convergences.
In both problems we construct exact solutions in terms of fundamental solutions of the Stokes equations \cite[Section IV.2]{Gal11}:
\begin{equation}\label{def:adjoint_deltas}
\mathbf{\Phi}(x):= 
\sum_{t\in\mathcal{E}}
\sum_{i=1}^{d}
\widetilde{\mathbf{T}}_{t}(x)\cdot \mathbf{e}_{i} ,\qquad 
\zeta(x):= 
\sum_{t\in\mathcal{E}}
\sum_{i=1}^{d}
\mathbf{T}_{t}(x)\cdot \mathbf{e}_{i},
\end{equation}
where $\{ \mathbf{e}_{i} \}_{i=1}^d$ denotes the canonical basis of $\mathbb{R}^{d}$ and 
\begin{align*}
\begin{array}{c}\displaystyle
\widetilde{\mathbf{T}}_{t}(\mathbf{x})
=
\left\{\begin{array}{ll}
-\dfrac{1}{4\pi}\bigg(\log|\mathbf{r}_{t}|\mathbb{I}^{}_2
-\dfrac{\mathbf{r}_{t}\mathbf{r}_{t}^{T}}{|\mathbf{r}_{t}|^2}
\bigg), &  d = 2, \\
\dfrac{1}{8\pi}\bigg(\dfrac{1}{|\mathbf{r}_{t}|}\mathbb{I}^{}_3
+\dfrac{\mathbf{r}_{t}\mathbf{r}_{t}^{T}}{|\mathbf{r}_{t}|^3}
\bigg), & d = 3;
\end{array}
\right. 
\displaystyle
\mathbf{T}_{t}(\mathbf{x})=\left\{\begin{array}{ll}
-\dfrac{\mathbf{r}_{t}}{2\pi|\mathbf{r}_{t}|^{2}}, & d = 2, \\
-\dfrac{\mathbf{r}_{t}}{4\pi|\mathbf{r}_{t}|^{3}}, & d = 3,
\end{array}
\right.
\end{array}
\end{align*} 
with $\mathbf{r}_{t} = x - t$. Here, $\mathbb{I}_d$ denotes the identity matrix in $\mathbb{R}^{d\times d}$.

Finally, we define $\mathbf{e}_{\mathbf{y}}:= \bar{\mathbf{y}}-\bar{\mathbf{y}}_h$, $e_p:=\bar{p}-\bar{p}_h$, $\mathbf{e}_{\mathbf{z}}:= \bar{\mathbf{z}}-\bar{\mathbf{z}}_h$, $e_r:=\bar{r}-\bar{r}_h$, $\mathbf{e}_{\mathbf{u}} := \bar{\mathbf{u}}-\bar{\mathbf{u}}_h$, and $\mathbf{e}_{\mathcal{U}} := \bar{\mathcal{U}}-\bar{\mathcal{U}}_h$.

\subsection{Optimization with Point Observations}\label{sec:num_ex_tracking}

The finite sequence of vectors 
$ 
\{ \mathbf{y}_{t} \}_{t \in \mathcal{D} } 
$ 
is computed from  the constructed solutions in such a way that the adjoint system \eqref{eq:adj_eq} holds. A straightforward computation reveals that, for $t\in\mathcal{D}$,
$
\mathbf{y}_{t}=\bar{\mathbf{y}}(t) - (\mathbf{e}_1 + \cdots + \mathbf{e}_d).
$
In order to simplify the construction of exact solutions, we have incorporated, in the \emph{momentum equation} of \eqref{eq:weak_pde}, an extra forcing term $\mathbf{f}\in\mathbf{L}^{\infty}(\Omega)$. With such a modification, the right hand side of the \emph{momentum equation} reads as follows: $(\mathbf{f}+\mathbf{u},\mathbf{v})_{\mathbf{L}^2(\Omega)}$. Finally, we will denote the total number of degrees of freedom as
$ 
\mathsf{Ndof}=2\dim(\mathbf{V}_{h})+2\dim(Q_{h})+\dim(\mathbf{U}_{h}).
$
\\~\\
\textbf{Example 1.} We let $\Omega = (0,1)^{2}$, $\alpha = 1.5 $, $\mathbf{a} = (-5,-5)^T$, $\mathbf{b} = (5,5)^T$, and
$
\mathcal{Z}=\{(0.25,0.25),(0.25,0.75),(0.75,0.25),(0.75,0.75)\}.
$
We define the exact optimal state as
\[\bar{\mathbf{y}}(x_1,x_2)=0.5\mathbf{curl}\left[(x_{1}x_{2}(1 - x_{1})(1 - x_{2}))^{2}\right],
\]
and $\bar{p}(x_{1},x_{2}) = x_{1}x_{2}(1-x_{1})(1-x_{2}) - {1}/{36}$,
while the exact optimal adjoint state is taken to be as in \eqref{def:adjoint_deltas}. 
\\~\\
\textbf{Example 2.} We set $\Omega=(0,1)^3$, $\mathbf{a} = (-10,-10,-10)^T$, $\mathbf{b} = (2,2,2)^T$, $\alpha = 1.99$, and $
\mathcal{Z}=\left\{(0.5,0.5,0.5)\right\}.
$
The exact optimal state is given by
\[
\bar{\mathbf{y}}(x_{1},x_{2},x_{3})  = -\frac{1}{\pi}\mathbf{curl}\big( (\sin(2\pi x_{1})\sin(2\pi x_{2})\sin(2\pi x_{3}) )^{2}\mathbf{e}_1 \big),
\]
and $\bar{p}(x_{1},x_{2},x_{3})  = x_{1}x_{2}x_{3} - {1}/{8}$.
The optimal adjoint state is as in \eqref{def:adjoint_deltas}.

\begin{figure}[ht]
\centering
\psfrag{error ad v L2}{{\huge $\|\mathbf{e}_{\mathbf{z}}\|_{\mathbf{L}^{2}(\Omega)}$}}
\psfrag{error st v Linf}{{\huge $\|\mathbf{e}_{\mathbf{y}}\|_{\mathbf{L}^{\infty}(\Omega)}$ }}
\psfrag{error st pr L2}{{\huge $\|e_{p}\|_{L^{2}(\Omega)}$ }}
\psfrag{error ct}{{\huge $\|\mathbf{e}_{\mathbf{u}}\|_{\mathbf{L}^{2}(\Omega)}$}}
\psfrag{error vel ad}{{\huge $\|\nabla \mathbf{e}_{\mathbf{z}}\|_{\mathbf{L}^{2}(\rho,\Omega)}$}}
\psfrag{Ndf-16}{\huge$\mathsf{Ndof}^{-1/6}$}
\psfrag{Ndf-10}{\huge$\mathsf{Ndof}^{-1}$}
\psfrag{Ndf-12}{\huge$\mathsf{Ndof}^{-1/2}$}
\psfrag{Ndf-23}{\huge$\mathsf{Ndof}^{-2/3}$}
\psfrag{Ndofs}{{\huge $\mathsf{Ndof}$}}
\begin{minipage}[c]{0.49\textwidth}
\centering
\psfrag{Example 1 - alfa 19 - errores}{\hspace{2.5cm}\huge{Convergence rates (Ex. 1)}}
\includegraphics[trim={0 0 0 0},clip,width=5.8cm,height=3.7cm,scale=0.58]{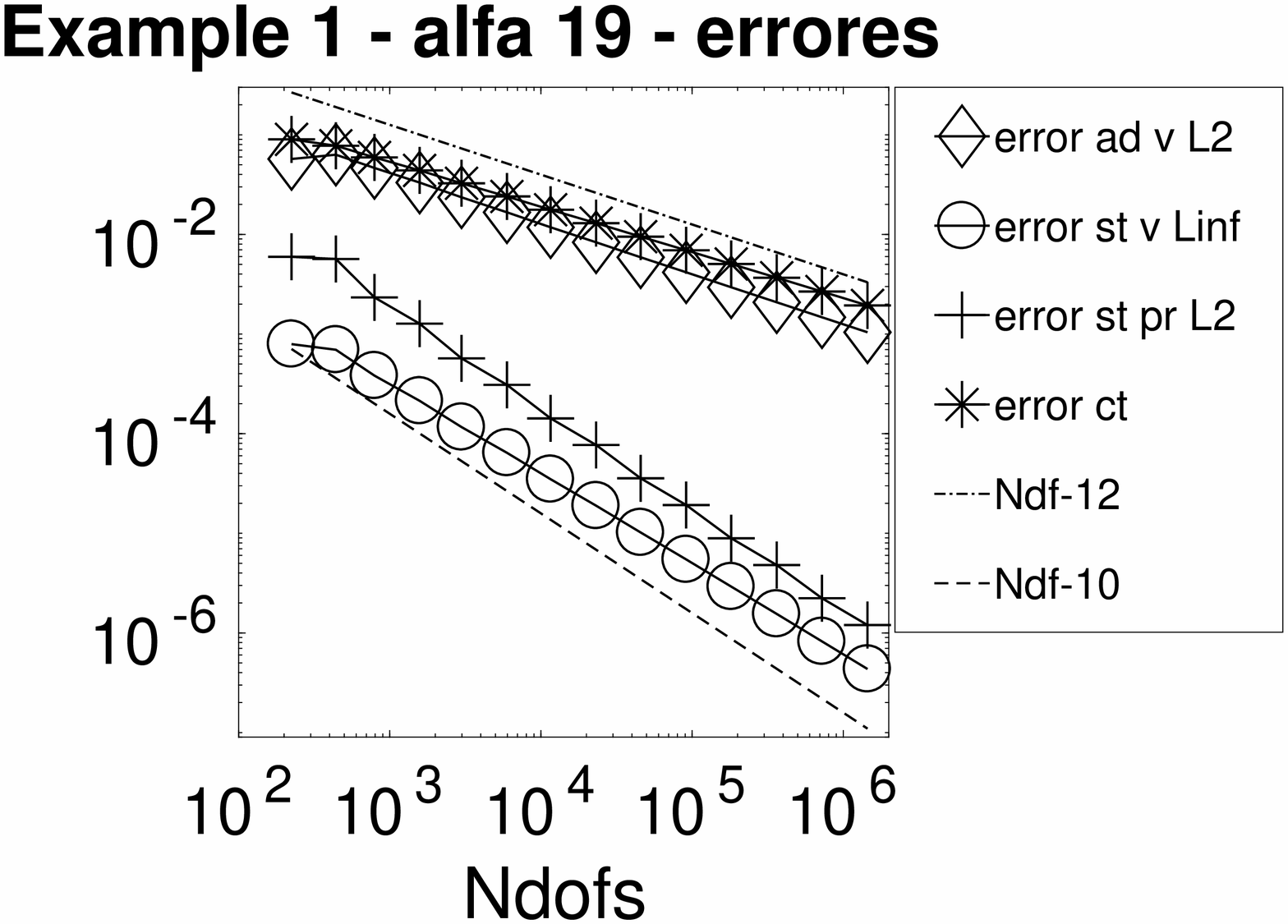} \\
\tiny{(A)}
\end{minipage}
\begin{minipage}[c]{0.49\textwidth}
\centering
\psfrag{Example 2 - alfa 199 - errores}{\hspace{2.8cm}\huge{Convergence rates (Ex. 2)}}
\includegraphics[trim={0 0 0 0},clip,width=5.8cm,height=3.7cm,scale=0.58]{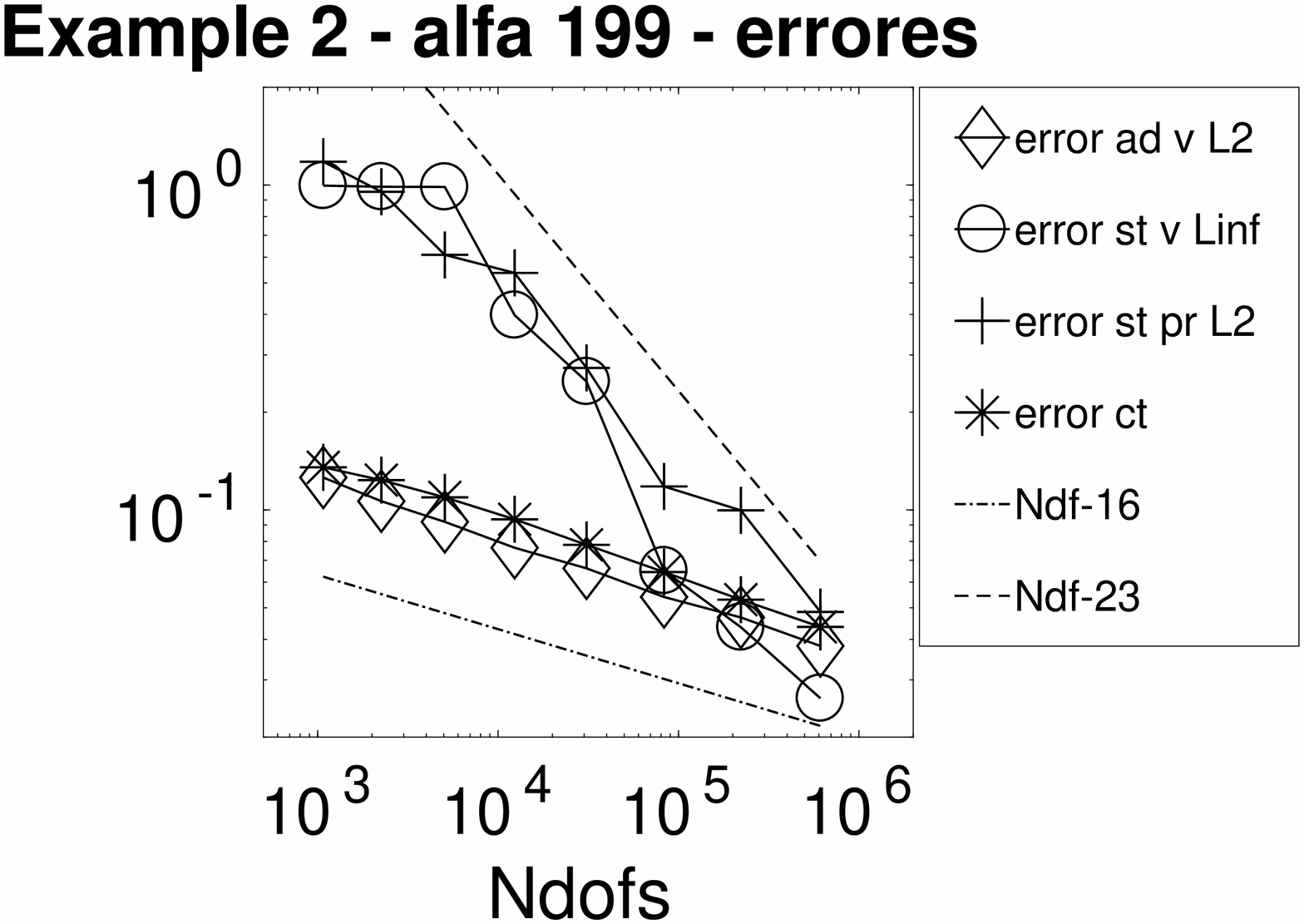} \\
\tiny{(B)}
\end{minipage}
\caption{Experimental rates of convergence for the approximation errors $\|\mathbf{e}_{\mathbf{z}}\|_{\mathbf{L}^{2}(\Omega)}$, $\|\mathbf{e}_{\mathbf{y}}\|_{\mathbf{L}^{\infty}(\Omega)}$, $\|e_{p}\|_{L^{2}(\Omega)}$, and $\|\mathbf{e}_{\mathbf{u}}\|_{\mathbf{L}^{2}(\Omega)}$ for Example 1 (A) and Example 2 (B).}
\label{fig:ex-1.1}
\end{figure}

We observe, in Fig. \ref{fig:ex-1.1}, that when approximating the optimal control variable $\bar{\mathbf{u}}$ and the adjoint velocity field $\bar{\mathbf{z}}$, the obtained experimental rates of convergence are in agreement with the estimates provided in \eqref{eq:global_reliability_2} and \eqref{global_reliability_adj_vel}, respectively. 

\subsection{Optimization with Singular Sources}\label{sec:num_ex_ss}

We now explore the performance of the discrete scheme \eqref{eq:discrete_state_eq_ss}--\eqref{eq:discrete_adj_eq_ss} with $d = 2$. In this case, the number of degrees of freedom is given by $ 
\mathsf{Ndof}=2\dim(\mathbf{V}_{h})+2\dim(Q_{h}) + 2l$, where $l = \#  \mathcal{D}$.
\\~\\
\textbf{Example 3.} We let $\Omega = (0,1)^{2}$, $\alpha = 1.99 $, $\mathbf{a}_{t} = (0,0)^T$, and $\mathbf{b}_{t} = (2,2)^T$ for all $t \in \mathcal{D}$ and
$
\mathcal{D}=\{(0.75,0.25)\}.
$
We define the exact optimal adjoint state as follows
\[\bar{\mathbf{z}}(x_1,x_2)=-\frac{4096}{27}\mathbf{curl}\left[(x_{1}x_{2}(1 - x_{1})(1 - x_{2}))^{2}\right],
\]
and $\bar{r}(x_{1},x_{2}) = x_{1}x_{2}(1-x_{1})(1-x_{2}) - {1}/{36}$. The exact optimal state is as in \eqref{def:adjoint_deltas}.

\begin{figure}[ht]
\centering
\psfrag{error st vel L2}{{\huge $\|\mathbf{e}_{\mathbf{y}}\|_{\mathbf{L}^{2}(\Omega)}$}}
\psfrag{error ad vel L2}{{\huge $\|\nabla\mathbf{e}_{\mathbf{z}}\|_{\mathbf{L}^{2}(\Omega)}$ }}
\psfrag{error ad pr L2}{{\huge $\|e_{r}\|_{L^{2}(\Omega)}$ }}
\psfrag{error ct}{{\huge $\VERT\mathbf{e}_{\mathcal{U}}\VERT_{\mathcal{D}}$}}
\psfrag{Ndf-16}{\huge$\mathsf{Ndof}^{-1/6}$}
\psfrag{Ndf-10}{\huge$\mathsf{Ndof}^{-1}$}
\psfrag{Ndf-12}{\huge$\mathsf{Ndof}^{-1/2}$}
\psfrag{Ndf-23}{\huge$\mathsf{Ndof}^{-2/3}$}
\psfrag{Ndofs}{{\huge $\mathsf{Ndof}$}}
\begin{minipage}[c]{0.47\textwidth}
\centering
\psfrag{Example 2 - alfa 199 - errores}{\hspace{2.5cm}\huge{Convergence rates (Ex. 3)}}
\includegraphics[trim={0 0 0 0},clip,width=6.0cm,height=3.7cm,scale=0.58]{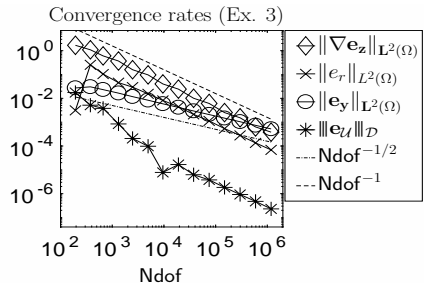} 
\end{minipage}
\caption{Experimental rates of convergence for the approximation errors $\|\nabla\mathbf{e}_{\mathbf{z}}\|_{\mathbf{L}^{2}(\Omega)}$, $\|e_{r}\|_{L^{2}(\Omega)}$, $\|\mathbf{e}_{\mathbf{y}}\|_{\mathbf{L}^{2}(\Omega)}$, and $\VERT\mathbf{e}_{\mathcal{U}}\VERT_{\mathcal{D}}$ within the setting of Example 3.}
\label{fig:ex-2}
\end{figure}

We observe, in Fig. \ref{fig:ex-2}, that when approximating the state velocity field $\bar{\mathbf{y}}$, the experimental rate of convergence for this variable is in agreement with the estimate provided in \EO{Theorem \ref{thm:velocity_rates_2}. We also observe that the experimental rate of convergence for the error approximation of the optimal control variable $\bar{\mathcal{U}}$ is in agreement with Theorem \ref{thm:control_rates_2}}. Finally we observe improved experimental rates of convergence for the adjoint variables $\bar{\mathbf{z}}$ and $\bar{r}$.

%
%

\footnotesize
\bibliographystyle{siam} 
\bibliography{biblio}

\end{document}